\documentclass[12pt]{article}

\usepackage{amssymb}
\usepackage{amsmath}
\usepackage{graphicx}
\usepackage{pstricks}
\usepackage{fullpage}
\usepackage{natbib}
\usepackage{jstyle}
\bibpunct{(}{)}{;}{a}{,}{,}
\renewcommand{\P}{\mbox{$\mathbb{P}$}}
\newcommand{\E}{\mbox{$\mathbb{E}$}}
\newcommand{\V}{\mbox{$\mathbb{V}$}}
\def\Var{{\mathbb V}\text{ar}}

\def\indicator{\mathbb I}
\def\u{{u}}
\def\H{{\mathcal H}}
\def\t{\top}
\def\A{{\mathcal A}}
\def\Amin{A_{\text{min}}}
\def\Amax{A_{\text{max}}}
\def\tr{\text{\,tr}\kern.03ex}
\newcommand{\ones}{\mbox{\boldmath$1$}}
\def\R{{\mathcal R}}
\newcommand{\LL}{h_{0}}
\def\overdot#1{\mbox{\it\.{#1}}} 
\def\gcv{\hat\R_{\text{gcv}}}
\def\GG{\mathcal G}
\def\asymplog{\approx}
\def\Dm{\nabla m}
\long\def\ignore#1{}
\def\Hh{H^{1/2}}
\def\vone{\ones}

\def\asymplog{\;\mathop{\mbox{\raise1.1ex\hbox{\tiny$\sim$}\kern-1.4ex\lower.1ex\hbox{\small
        $\asymp$}}}\;}
\let\hat\widehat
\let\tilde\widetilde

\begin{document}

\begin{center}
\begin{tabular}{c}
\\[.5in]
{\LARGE\bf Rodeo: Sparse Nonparametric Regression} \\[5pt]
{\LARGE\bf in High Dimensions}\\[20pt]
\end{tabular}
\begin{tabular}{ccc}
\large John Lafferty\footnotemark  && \large Larry Wasserman\footnotemark \\[2pt]
Computer Science Department && Department of Statistics \\[-1pt]
Machine Learning Department && Machine Learning Department \\[-1pt]
Carnegie Mellon University && Carnegie Mellon University \\[10pt]
\end{tabular}

\today\\[10pt]
\end{center}

\vskip30pt
\centerline{\bf Abstract}
\begin{quote}
  
We present a greedy method for simultaneously performing local
bandwidth selection and variable selection in nonparametric
regression.  The method starts with a local linear estimator with
large bandwidths, and incrementally decreases the bandwidth of
variables for which the gradient of the estimator with respect to
bandwidth is large.  The method---called {\it rodeo\/} (regularization
of derivative expectation operator)---conducts a sequence of
hypothesis tests to threshold derivatives, and is easy to implement.
Under certain assumptions on the regression function and sampling
density, it is shown that the rodeo applied to local linear smoothing
avoids the curse of dimensionality, achieving near optimal minimax
rates of convergence in the number of relevant variables, as if these
variables were isolated in advance.

\end{quote}


\vskip10pt

\footnotetext[1]{Research supported by NSF grants 
CCR-0122481, IIS-0312814, and IIS-0427206.}
\footnotetext[2]{Research supported by 
NIH grants R01-CA54852-07 and MH57881 and NSF grant DMS-0104016.}

\clearpage

\section{Introduction}

Estimating a high dimensional regression function is notoriously
difficult due to the curse of dimensionality.  Minimax theory
precisely characterizes the curse.  Let 
\begin{equation}
Y_i = m(X_i) + \epsilon_i, \ \ \ i=1, \ldots, n
\end{equation}
where $X_i = (X_{i}(1), \ldots, X_{i}(d))\in
\mathbb{R}^d$ is a $d$-dimensional covariate, $m: \mathbb{R}^d \to
\mathbb{R}$ is the unknown function to estimate, and
$\epsilon_i \sim N(0,\sigma^2)$.
Then if $m$ is in $W_2(c)$, the $d$-dimensional Sobolev ball of order two
and radius $c$, it is well known that
\begin{equation}
\liminf_{n\to\infty} n^{4/(4+d)}  \;
\inf_{\widehat{m}_n} \; \sup_{m\in W_2(c)} \; \R(\hat{m}_n, m) \;>\; 0
\,,
\end{equation}
where $\R(\widehat{m}_n, m) = \E_m\int (\widehat{m}_n(x) -
m(x))^2 \,dx$ is the risk of the estimate $\widehat{m}_n$ constructed
\relax from a sample of size $n$ \citep{gyorfi:2002,Ston:Hans:Koop:Truo:1997}.
Thus, the best rate of convergence is $n^{-4/(4+d)}$, which is
impractically slow if $d$ is large.  

However, for some applications it is reasonable to expect that the
true function only depends on a small number of the total
covariates.   Suppose that
$m$ satisfies such a sparseness condition, so that
\begin{equation}\label{eq::sparse}
m(x) = m(x_R) 
\end{equation}
where $x_R=(x_j:\ j\in R)$, $R\subset \{1, \ldots, d\}$ is a subset of
the $d$ covariates, of size $r=|R| \ll d$.  We call $\{x_j\}_{j\in R}$
the {\it relevant variables\/}.  Note that if an oracle were to
identify and isolate the relevant variables, the better minimax rate
of $n^{-4/(4+r)}$ could be achieved, and this would be the fastest
rate possible.  Thus, we are faced with the problem of variable
selection in nonparametric regression.  Our strategy is to seek a
greedy method that incrementally searches through bandwidths in small steps.

A large body of previous work has addressed this fundamental problem,
which has led to a variety of methods to combat the curse of
dimensionality. Many of these are based on very clever, though often
heuristic techniques.  For additive models of the form $m(x) = \sum_j
m_j(x_j)$, standard methods like stepwise selection, $C_p$ and AIC can
be used \citep{Hastie:2001}.  For spline models, \cite{Zhang:04} use
likelihood basis pursuit, essentially the lasso adapted to the spline
setting.  CART \citep{cart:84} and MARS \citep{Friedman:91}
effectively perform variable selection as part of their function
fitting.  Support vector regression can be seen as creating a sparse
representation using basis pursuit in a reproducing kernel Hilbert
space \citep{Girosi:97}. There is also a large literature on Bayesian
methods, including methods
for sparse Gaussian processes \citep{Tipping:01,Smola:01,Lawrence:03};
see \cite{George:97} for a brief survey.
More recently, \cite{Li:05} use independence testing for variable
selection and \cite{buhl:2005} introduced a boosting approach.  While
these methods have met with varying degrees of empirical success, they
can be challenging to implement and demanding computationally.
Moreover, these methods are typically very difficult to analyze
theoretically, and so come with limited formal guarantees.  Indeed, the
theoretical analysis of sparse {\it parametric\/} estimators such as
the lasso \citep{Tibshirani:96} is challenging, and only recently has
significant progress been made on this front in the
statistics and signal processing communities
\citep{Donoho:04,Knight:00,Tropp:04,Tropp:06,Fan:Peng:2004,Fan:Li:2001}.

In this paper we present a new approach for sparse nonparametric
function estimation that is both computationally simple and amenable
to theoretical analysis.  We call the general framework {\it rodeo\/},
for ``regularization of derivative expectation operator.''  It is based on
the idea that bandwidth and variable selection can be simultaneously
performed by computing the infinitesimal change in a nonparametric
estimator as a function of the smoothing parameters, and then
thresholding these derivatives to get a sparse estimate.
As a simple version of this principle we use hard thresholding,
effectively carrying out a sequence of hypothesis tests.  A modified
version that replaces testing with soft thresholding may be viewed as
solving a sequence of lasso problems.  The potential appeal of this
approach is that it can be based on relatively simple and
theoretically well understood nonparametric techniques such as local
linear smoothing, leading to methods that are simple to implement and
can be used in high dimensional problems. Moreover, we show that they
can achieve near optimal minimax rates of convergence, and therefore
circumvent the curse of dimensionality when the true function is
indeed sparse.  When applied in one dimension, our method yields a
local bandwidth selector and
is similar to the estimators of \cite{Rupp:1997} and \cite{Lepski:1997}.
The method in \cite{Lepski:1997} and its multivariate extension in 
Kerkyacharian, Lepski and Picard (2001) 
yield estimators that are 
more refined than our method in the sense that their estimator is spatially
adaptive over large classes of function spaces.
However, their method is not greedy: it involves searching over
a large class of bandwidths. Our goal is to develop a greedy method
that scales to high dimensions.

Our method is related to the structural adapation method of
\cite{Hristache:01} and \cite{Samarov:05}, which is designed for
multi-index models.  The general multi-index model is
\begin{equation}\label{eq::index}
Y=g_0 (Tx) + \epsilon
\end{equation}
where $x\in\mathbb{R}^d$ and
$T$ is a linear orthonormal mapping from
$\mathbb{R}^d$ onto $\mathbb{R}^r$ with $r < d$.
Variable selection corresponds to taking $T$
to be a $r$ by $d$ matrix of 0's and 1's with each
$T_{ij}=1$ if $x_j$ is the $i^{\rm th}$ relevant variable.
Nonparametric variable selection can also be regarded
as a special case of the partially linear model in
\cite{Samarov:05}, which takes
\begin{equation}\label{eq::partlin}
Y= \theta^T x_1 + G(x_2) + \epsilon
\end{equation}
where $x=(x_1,x_2)$.  Taking $\theta$ to be zero yields the model in
this paper.  The advantage of structural adaptation is that it yields,
under certain conditions, $\sqrt{n}$ estimates of the image of $T$ in
(\ref{eq::index}) and $\theta$ in (\ref{eq::partlin}).  However,
structural adaptation does not yield optimal bandwidths or optimal
estimates of the regression function, although this 
is not the intended goal of the method.

In the following section we outline the basic rodeo approach, which is
actually a general strategy that can be applied to a wide range of
nonparametric estimators.  We then specialize in
Section~\ref{sec::multivariate} to the case of local linear smoothing,
since the asymptotic properties of this smoothing technique are fairly
well understood.  In particular, we build upon the analysis of
\cite{Ruppert:1994} for local linear regression; a notable difference
is that we allow the dimension to increase with sample size,
which requires a more detailed analysis of the asymptotics.  In
Section~\ref{sec::examples} we present some simple examples of the
rodeo, before proceeding to an analysis of its properties in
Section~\ref{sec::properties}.  Our main theoretical result
characterizes the asymptotic running time, selected bandwidths, and
risk of the algorithm.  
Finally, in Section~\ref{sec::mods} we present
further examples and discuss several extensions of the basic version
of the rodeo considered in the earlier sections.
The proofs of the theoretical properties of the rodeo are
given in Section~\ref{sec::proofs}.

\section{Rodeo: The Main Idea}

The key idea in our approach is as follows.
Fix a point $x$ and
let $\hat{m}_h(x)$ denote an estimator
of $m(x)$ based on a vector
of smoothing parameters
$h=(h_1, \ldots, h_d)$.
If $c$ is a scalar,
then we write
$h=c$ to mean
$h=(c, \ldots, c)$.

Let $M(h)=\E(\hat{m}_h(x))$ denote the mean
of $\hat{m}_h(x)$.
For now, assume that $x=x_i$ is one of the observed data points
and that $\hat{m}_0(x) = Y_i$.
In that case,
$m(x) = M(0)=\E(Y_i)$.
If $P=(h(t):\  0 \leq t \leq 1)$ is a smooth path through the set of
smoothing parameters 
with $h(0) = 0$ and $h(1) = 1$
(or any other fixed, large bandwidth)
then
\begin{subequations}
\begin{eqnarray}
m(x) &=& M(0) = M(1) + M(0) - M(1) \\
&=&  M(1) - \int_0^1 \frac{dM(h(s))}{ds} \, ds\\
\label{eq::basic}
&=&  M(1) - \int_0^1 \bigl\langle D(s), \overdot{h}(s)\bigr\rangle ds
\end{eqnarray}
\end{subequations}
where
\begin{equation}
D(h) = \nabla M(h) = 
\left(\frac{\partial M}{\partial h_1},\ldots ,
     \frac{\partial M}{\partial h_d} \right)^T
\end{equation}
is the gradient of $M(h)$ and
$\overdot{h}(s) = \frac{dh(s)}{ds}$ is the derivative of $h(s)$ along the path.
A biased, low variance estimator of $m(x)$ is $\hat{m}_1(x)$.
An unbiased estimator of $D(h)$ is
\begin{equation}\label{eq::Z}
Z(h) = \left(\frac{\partial\hat{m}_h(x)}{\partial h_1}, \ldots,
              \frac{\partial\hat{m}_h(x)}{\partial h_d} \right)^T.
\end{equation}
The naive estimator
\begin{equation}
\hat{m}(x)=\hat{m}_1(x) -  \int_0^1 \bigl\langle Z(s), \overdot{h}(s)\bigr\rangle ds
\end{equation}
is identically equal to $\hat{m}_0(x)=Y_i$,
which has poor risk since
the variance of $Z(h)$ is large for small $h$.
However, our sparsity assumption on $m$ suggests that there should be paths
for which $D(h)$ is also sparse.
Along such a path, we replace $Z(h)$
with an estimator $\hat{D}(h)$ that makes use of the sparsity assumption.
Our estimate of $m(x)$ is then
\begin{equation}
\tilde{m}(x)  =  \hat{m}_1(x) - \int_0^1 \bigl\langle \hat{D}(s), \overdot{h}(s)\bigr\rangle ds\,.
\end{equation}
To implement this idea we need to do two things:
(i) we need to find a path for which the derivative is sparse and (ii) we need to
take advantage of this sparseness when estimating $D$ along that
path.

\begin{figure}[t]
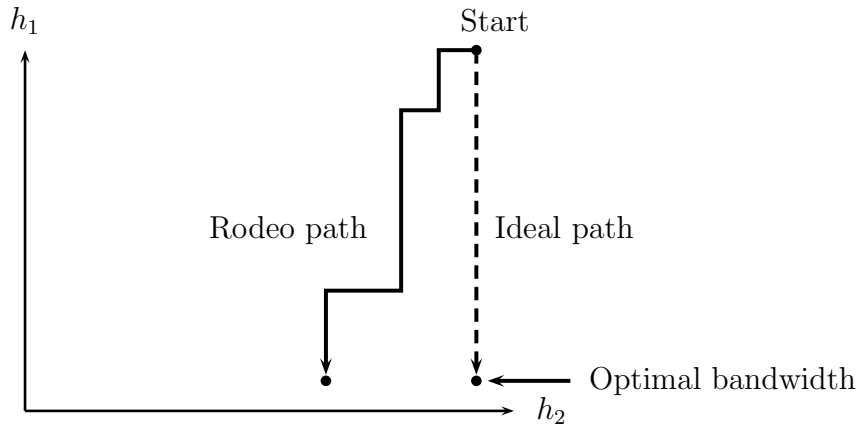

\vskip5pt
\psset{xunit=.5cm,yunit=.4cm}
\pspicture(-6,0)(14,16)
\psline[linewidth=1pt]{->}(2,2)(2,14)
\psline[linewidth=1pt]{->}(2,2)(15,2)
\psline[linewidth=1.5pt,linestyle=dashed]{->}(14,14)(14,3.2)
\psline[linewidth=1.5pt]{->}(14,14)(13,14)(13,12)(12,12)(12,11)(12,6)(11,6)(10,6)(10,3.2)
\psdot[dotsize=4pt](14,3)
\psdot[dotsize=4pt](14,14)
\psdot[dotsize=4pt](10,3)
\rput(16,2){$h_2$}
\rput(2,15){$h_1$}
\rput(14.5,15.0){Start}
\rput[l](17,3){Optimal bandwidth}
\psline[linewidth=1.5pt]{->}(16.5,3)(14.3,3)
\rput[l](14.5,8){Ideal path}
\rput(9,8){Rodeo path}
\endpspicture
\vskip-10pt
\caption{Conceptual illustration: The bandwidths for the relevant
variables ($h_1$) are shrunk, while the bandwidths for the irrelevant
variables ($h_2$) are kept relatively large.}
\vskip10pt
\label{fig::artwork}
\end{figure}

The key observation is that if $x_j$ is irrelevant, then we expect
that changing the bandwidth $h_j$ for that variable should cause only
a small change in the estimator $\hat{m}_h(x)$.  Conversely, if $x_j$
is relevant, then we expect that changing the bandwidth $h_j$ for that
variable should cause a large change in the estimator.  Thus, $Z_j =
{\partial \widehat{m}_h(x)}/{\partial h_j}$ should discriminate
between relevant and irrelevant covariates.  To simplify the
procedure, we can replace the continuum of bandwidths with a
discrete set where each $h_j \in {\cal B} = \{\LL, \beta\LL, \beta^2\LL, \ldots \}$ for some $0 <
\beta < 1$.  Moreover, we can proceed in a greedy fashion by
estimating $D(h)$ sequentially with $h_j \in {\cal B}$ and setting
$\hat{D}_j(h)=0$ when $h_j < \hat{h}_j$, where $\hat{h}_j$ is the first $h$
such that $|Z_j(h)| < \lambda_j(h)$ for some threshold $\lambda_j$.  
This greedy version, coupled with
the hard threshold estimator, yields $\tilde{m}(x) =
\hat{m}_{\hat{h}}(x)$.  A conceptual illustration of the idea is shown
in Figure~\ref{fig::artwork}.

To further elucidate the idea,
consider now the one-dimensional case $x\in\mathbb{R}$, so that
\begin{equation}
m(x) = M(1) - \int_0^1 \frac{ dM(h)}{dh} dh = M(1) - \int_0^1 D(h)\, dh.
\end{equation}
Suppose that
$\hat{m}_h(x)=\sum_{i=1}^n Y_i \,\ell_i(x,h)$ 
is a linear estimator, where the weights
$\ell_i(x,h)$ depend on a bandwidth $h$.

In this case
\begin{equation}
Z(h) = \sum_{i=1}^n Y_i\, \ell_i'(x,h)
\end{equation}
where the prime denotes differentiation with respect to $h$.
Then we set
\begin{equation}
\tilde{m}(x) = \hat{m}_1(x)  - \int_0^1 \hat{D}(h) \, dh
\end{equation}
where $\hat{D}(h)$ is an estimator of $D(h)$.
Now, 
\begin{equation}
Z(h) \approx N(b(h),s^2(h))
\end{equation}
where, for typical smoothers,
$b(h) \approx Ah$ and
$s^2(h) \approx C/nh^3$
for some constants $A$ and $C$.
Take the hard threshold estimator
\begin{equation}
\hat{D}(h) = Z(h) I\bigl(|Z(h)| > \lambda(h)\bigr)
\end{equation}
where $\lambda(h)$ is chosen to be slightly larger than
$s(h)$.
An alternative is the soft-threshold estimator
\begin{equation}
\hat{D}(h) = {\rm sign}(Z(h)) (|Z(h)| - \lambda(h))_+ .
\end{equation}

The greedy algorithm, coupled with the hard threshold estimator, yields
a bandwidth selection procedure based on testing.
This approach to bandwidth selection
is very similar to that 
of \cite{Lepski:1997}, who take
\begin{equation}
\hat{h} = \max\{ h\in {\cal H}:\ \phi(h,\eta)=0\ \mbox{for all}\ \eta < h\}
\end{equation}
where $\phi(h,\eta)$ is a test for whether $\hat{m}_\eta$ improves
on $\hat{m}_h$.  This more refined test leads to estimators that
achieve good spatial adaptation over large function classes.  
Kerkyacharian, Lepski and Picard (2001) extend the idea to multiple dimensions.
Our approach is also similar to a method of \cite{Rupp:1997} that uses a
sequence of decreasing bandwidths and then estimates the optimal
bandwidth by estimating the mean squared error as a function of
bandwidth. Our greedy approach only tests whether an infinitesimal
change in the bandwidth from its current setting leads to a
significant change in the estimate, and is more easily extended to a
practical method in higher dimensions.

\section{Rodeo Using Local Linear Regression}
\label{sec::multivariate}

Now we present the multivariate rodeo in detail.
We use local linear smoothing as the basic method since
it is known to have many good properties.
Let $x = (x(1),\ldots, x(d))$ be some target point at which we want to
estimate $m$.
Let $\hat{m}_H(x)$ denote the local linear estimator of
$m(x)$ using bandwidth matrix $H$.
Thus,
\begin{equation}
\hat{m}_H(x) = e_1^T (X_x^T W_x X_x)^{-1} X_x^T W_x Y \equiv S_x Y
\end{equation}
where
$e_1 = (1,0,\ldots, 0)^T$,
\begin{equation}
X_x = \left(\begin{array}{cc} 1 & (X_1-x)^T \\ \vdots & \vdots \\ 1 & (X_n-x)^T\end{array}\right),
\end{equation}
$W_x$ is diagonal with
$(i,i)$ element
$K_H(X_i-x)$ and
$K_H(u) = |H|^{-1/2} K(H^{-1/2} u)$.
The estimator
$\hat{m}_H$ can be written as
\begin{equation}
\hat{m}_H(x) = \sum_{i=1}^n G(X_i,x,h) \, Y_i
\end{equation}
where
\begin{equation}
G(u,x,h) = e_1^T (X_x^T W_x X_x)^{-1} 
\left(\begin{array}{cc}1 \\ (u-x)^T \end{array}\right) K_H(u-x)
\end{equation}
is called the \emph{effective kernel}.
One can regard local linear regression
as a refinement of kernel regression
where the effective kernel $G$ adjusts for boundary bias and design bias;
see \cite{Fan:1992}, \cite{Hast:Load:1993} and \cite{Ruppert:1994}.

We assume that the covariates are random with density $f(x)$
and that $x$ is interior to the support of $f$.
We make the same assumptions as Ruppert and Wand (1994) in their
analysis of the bias and variance of local linear regression. In
particular:
\begin{enumerate}
\item[(i)] The kernel $K$ has compact support with zero odd
moments and there exists $\nu_2 =\nu_2(K) \neq 0$ such that
\begin{equation}
\label{eq:nudef}
\int \u \u^\t K(\u)\,d\u = \nu_2(K) I
\end{equation}
where $I$ is the $d\times d$ identity matrix.
\item[(ii)] The sampling density $f(x)$ is continuously
differentiable and strictly positive.
\end{enumerate}

In the version of the algorithm that follows,
we take $K$ to be a product kernel and
$H$ to be diagonal with elements
$h=(h_1, \ldots, h_d)$ and we write
$\hat{m}_h$ instead of $\hat{m}_H$.

Our method is based on the statistic
\begin{equation}\label{eq::Zj}
Z_j = \displaystyle \frac{\partial \widehat{m}_h(x)}{\partial h_j} = \sum_{i=1}^n G_j(X_i,x,h) Y_i
\end{equation}
where
\begin{equation}
G_j(u,x,h) = \frac{\partial G(u,x,h)}{\partial h_j}.
\end{equation}
Let
\begin{equation}
\mu_j \equiv \mu_j(h) = \E(Z_j|X_1, \ldots, X_n) = \sum_{i=1}^n G_j(X_i,x,h) m(X_i)
\end{equation}
and
\begin{equation}\label{eq::sj}
s_j^2 \equiv s_j^2 (h)= \V(Z_j|X_1, \ldots, X_n) = \sigma^2 \sum_{i=1}^n G_j(X_i,x,h)^2 .
\end{equation}
In Section \ref{sec::sigma} we explain how to estimate $\sigma$; for
now, assume that $\sigma$ is known.  The hard thresholding version of
the rodeo algorithm is described in Figure \ref{fig::hard-rodeo}.
We make use of a sequence $c_n$ satisfying 
$d c_n = \Omega(\log n)$, where we write $f(n) = \Omega(g(n))$ if 
$\liminf_{n} \left|\frac{f(n)}{g(n)}\right| > 0$.

\begin{figure*}[ht]
{\it Rodeo: Hard thresholding version}
\vskip5pt
\hrule
\vskip5pt
\begin{enumerate}
\item {\it Select\/} 
$\beta_n = n^{-\alpha/\log^3 n}$
for some $0 < \alpha < 1$
and initial bandwidth
\begin{equation}
\label{eq:initial_h}
\LL = \frac{c_0}{\log \log n}
\end{equation}
for some constant $c_0$.  
Let $c_n$ be a sequence satisfying $c_n = O(1)$.
\item {\it Initialize\/} the bandwidths, and activate all covariates:
\begin{enumerate}
\item $h_j = \LL$, $j=1,2,\ldots, d$.
\item $\A = \{1,2,\ldots, d\}$
\end{enumerate}
\item {\it While $\A$ is nonempty}, do for each $j\in \A$:
\begin{enumerate}
\item Compute the estimated derivative expectation: $Z_j$ (equation \ref{eq::Zj})
and $s_j$ (equation \ref{eq::sj}).
\item Compute the threshold
$ \lambda_j = \displaystyle s_j \sqrt{2 \log (n c_n)}$.
\item If $|Z_j| > \lambda_j$, then set $h_j \leftarrow \beta h_j$;
  otherwise remove $j$ from $\A$.
\end{enumerate}
\item {\it Output} bandwidths $h^\star=(h_1,\ldots,h_d)$ and estimator
$\tilde{m}(x) = \hat{m}_{h^\star}(x)$.
\end{enumerate}
\vskip10pt
\hrule
\vskip10pt
\caption{The hard thresholding version of the rodeo, which can be applied
using the derivatives $Z_j$ of any nonparametric smoother.}
\label{fig::hard-rodeo}
\vskip15pt
\end{figure*}

To derive an explicit expression for $Z_j$, equivalently $G_j$, we use
\begin{equation}
\frac{\partial A^{-1}}{\partial h} = - A^{-1} \frac{\partial A}{\partial h} A^{-1}
\end{equation}
to get that 
\begin{subequations}
\begin{eqnarray}
Z_j &=& \frac{\partial\hat m_h(x)}{\partial h_j} \\
&=&
\nonumber
e_1^\t  (X^\t W X)^{-1} X^\t \frac{\partial W}{\partial h_j} Y - 
e_1^\t  (X^\t W X)^{-1} X^\t \frac{\partial W}{\partial h_j} X (X^\t W
X)^{-1} X^\t W Y \\
&& \\
&=& e_1^\t  (X^\t W X)^{-1} X^\t \frac{\partial W}{\partial h_j} (Y -
X\hat \alpha)
\end{eqnarray}
\end{subequations}
where $\hat\alpha = (X^\t W X)^{-1} X^\t W Y$ is the coefficient
vector for the local linear fit (and we have dropped the dependence on
the local point $x$ in the notation).

Note that the factor $|H|^{-1} = \prod_{i=1}^d 1/h_i$ in the
kernel cancels in the expression for $\hat m$, and therefore we can
ignore it in our calculation of $Z_j$.  Assuming a product kernel we have
\begin{equation}
W \;=\; \text{diag}\left(
\prod_{j=1}^d K((X_{1j}-x_j)/h_j),\ldots,\prod_{j=1}^d K((X_{nj}-x_j)/h_j)\right)
\end{equation}
and ${\partial W}/{\partial h_j} = W L_j$ where 
\begin{equation}
L_j = \text{diag}\left(
\frac{\partial \log K((X_{1j}-x_j)/h_j)}{\partial h_j}, \ldots, 
\frac{\partial \log K((X_{nj}-x_j)/h_j)}{\partial h_j}\right)
\end{equation}
and thus
\begin{eqnarray}
\label{eq:inverse}
Z_j 
&=& e_1^\t  (X^\t W X)^{-1} X^\t W L_j (Y - X\hat \alpha) =
e_1^\t B L_j (I - X B) Y = G_j(x,h)^\t Y
\end{eqnarray}
where $B = (X^\t W X)^{-1} X^\t W$.

The calculation of $L_j$ is typically 
straightforward.  As two examples, with the Gaussian kernel $K(u) = \exp(-u^2/2)$ we have
\begin{equation}
L_j = \frac{1}{h_j^3}\, \text{diag}\left(
(X_{1j}-x_j)^2, \ldots, (X_{nj}-x_j)^2\right)
\end{equation}
and for the Epanechnikov kernel $K(u) = (5-x^2)\, \indicator(|x| \leq
\sqrt{5})$ we have
\begin{subequations}
\begin{eqnarray}
\lefteqn{L_j \;=\;
\frac{1}{h_j^3}\, \text{diag}\left(
\frac{2(X_{1j}-x_j)^2}{5-(X_{1j}-x_j)^2/h_j^2} 
\, \indicator(|X_{1j}-x_j| \leq \sqrt{5} h_j), 
\ldots, \right.} && \\
&& \hskip60pt
\left.\frac{2(X_{nj}-x_j)^2}{5-(X_{nj}-x_j)^2/h_j^2} 
\, \indicator(|X_{1j}-x_j| \leq \sqrt{5} h_j)\right)
\end{eqnarray}
\end{subequations}

\section{Examples}
\label{sec::examples}

In this section we illustrate the rodeo on some examples.  We return
to the examples later when we discuss estimating $\sigma$, as well as
a global (non-local) version of the rodeo.

\subsection{Two Relevant Variables}
\label{sec::tworel}

\begin{figure*}
\begin{center}
\begin{tabular}{ccc}
&& \\[-30pt]
\includegraphics[width=2.0in,angle=0]{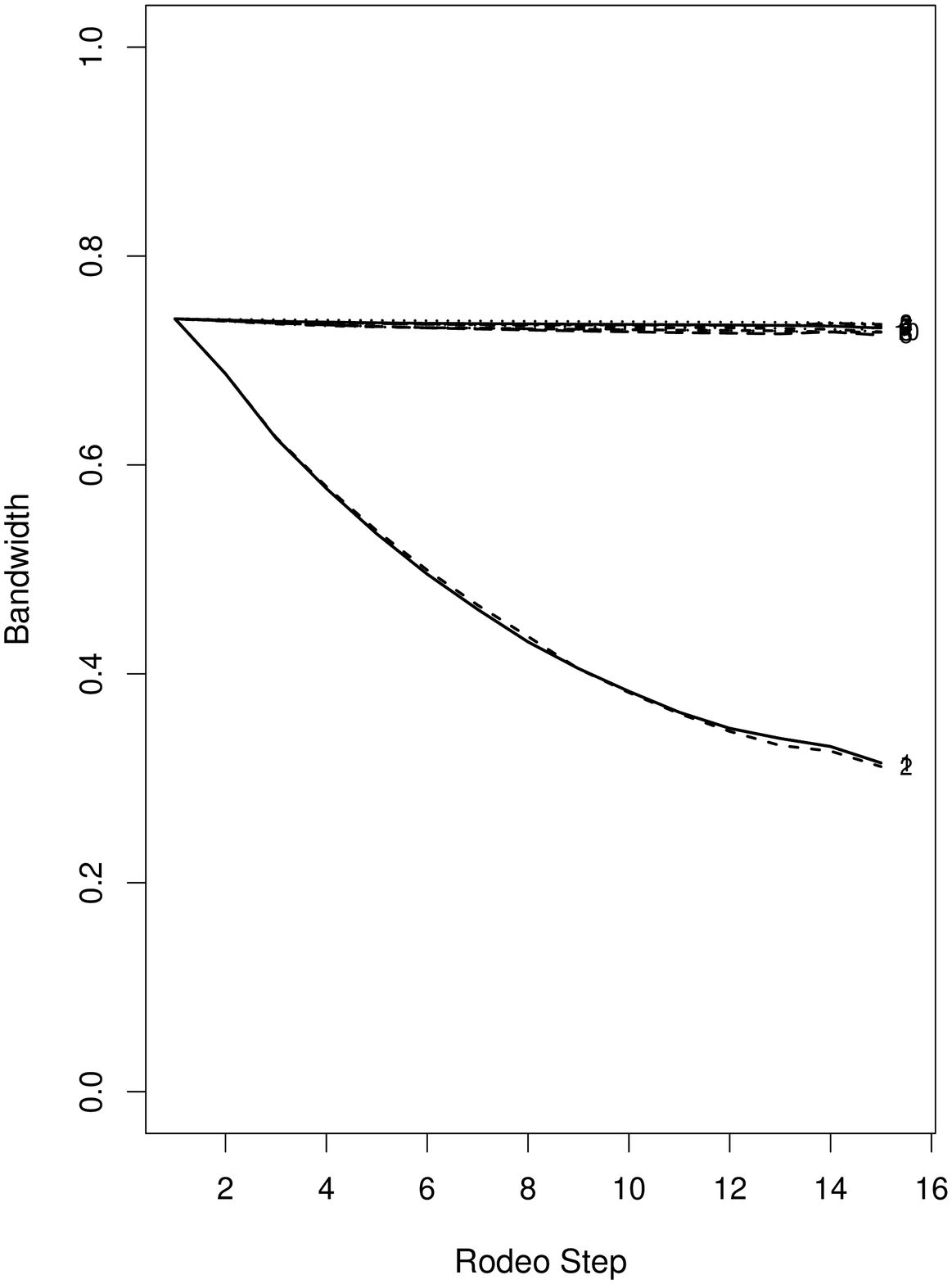} &
\includegraphics[width=2.0in,angle=0]{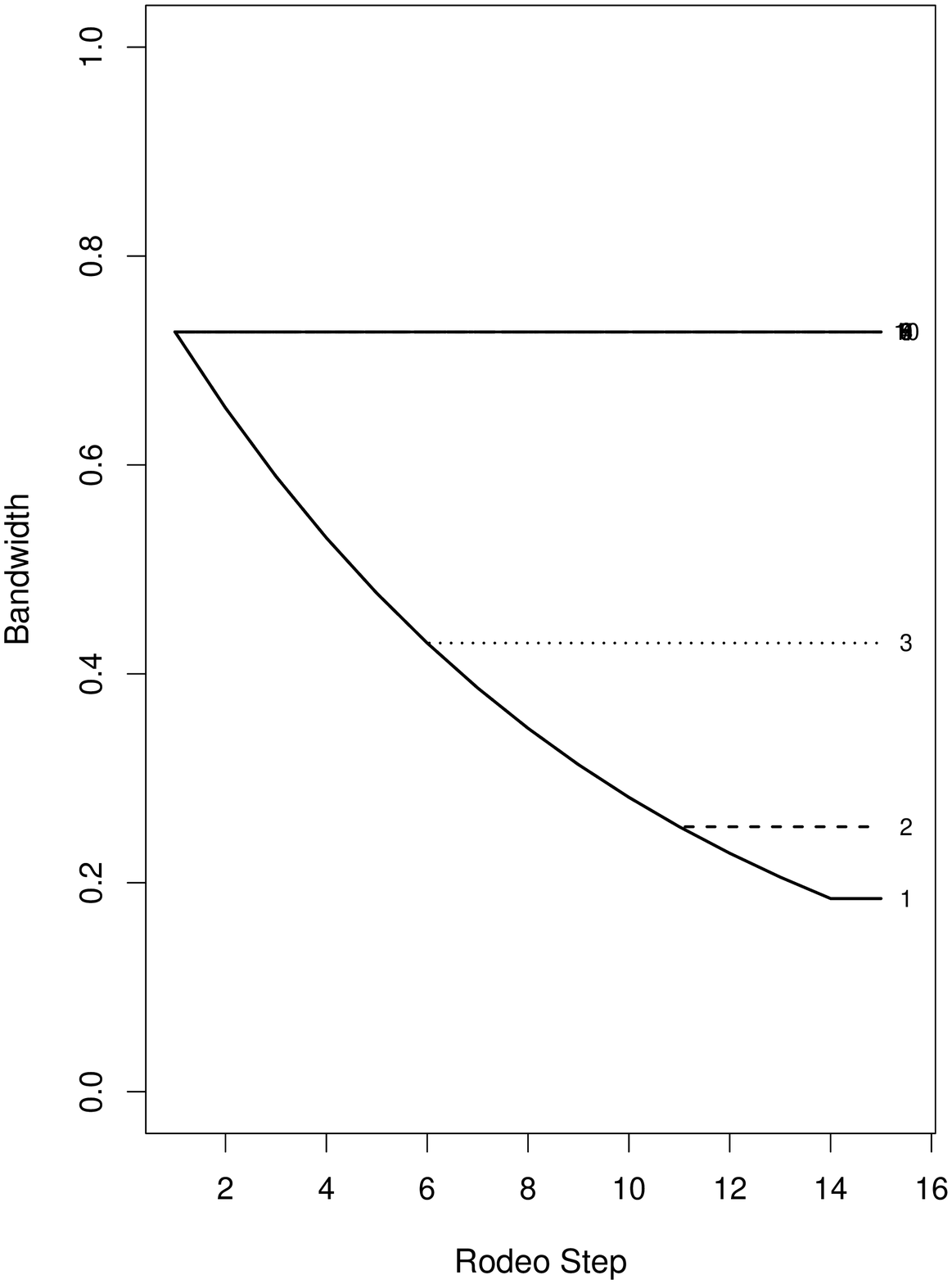} &
\includegraphics[width=2.0in,angle=0]{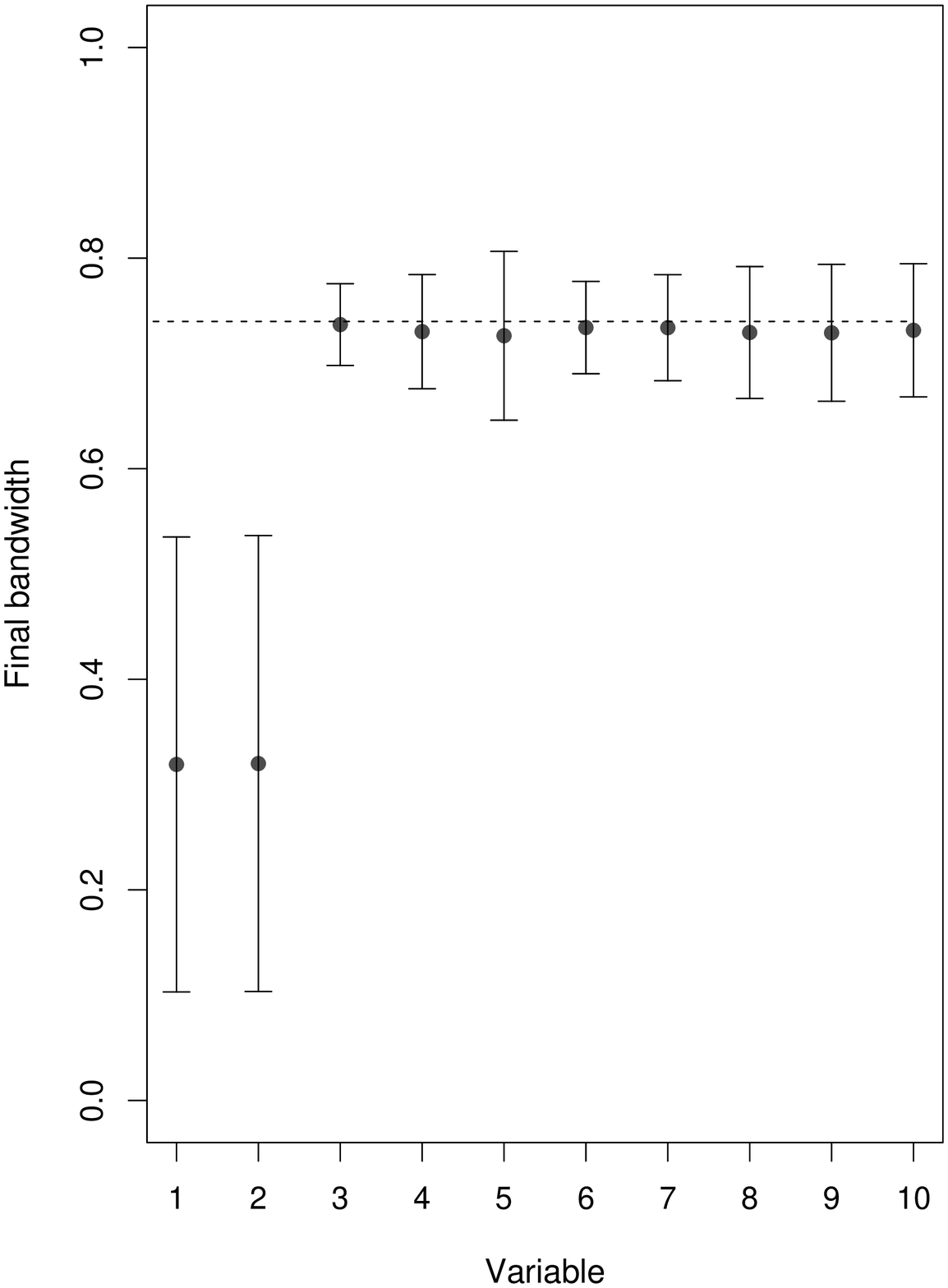} \\
&& \\[-30pt]
\includegraphics[width=2.0in,angle=0]{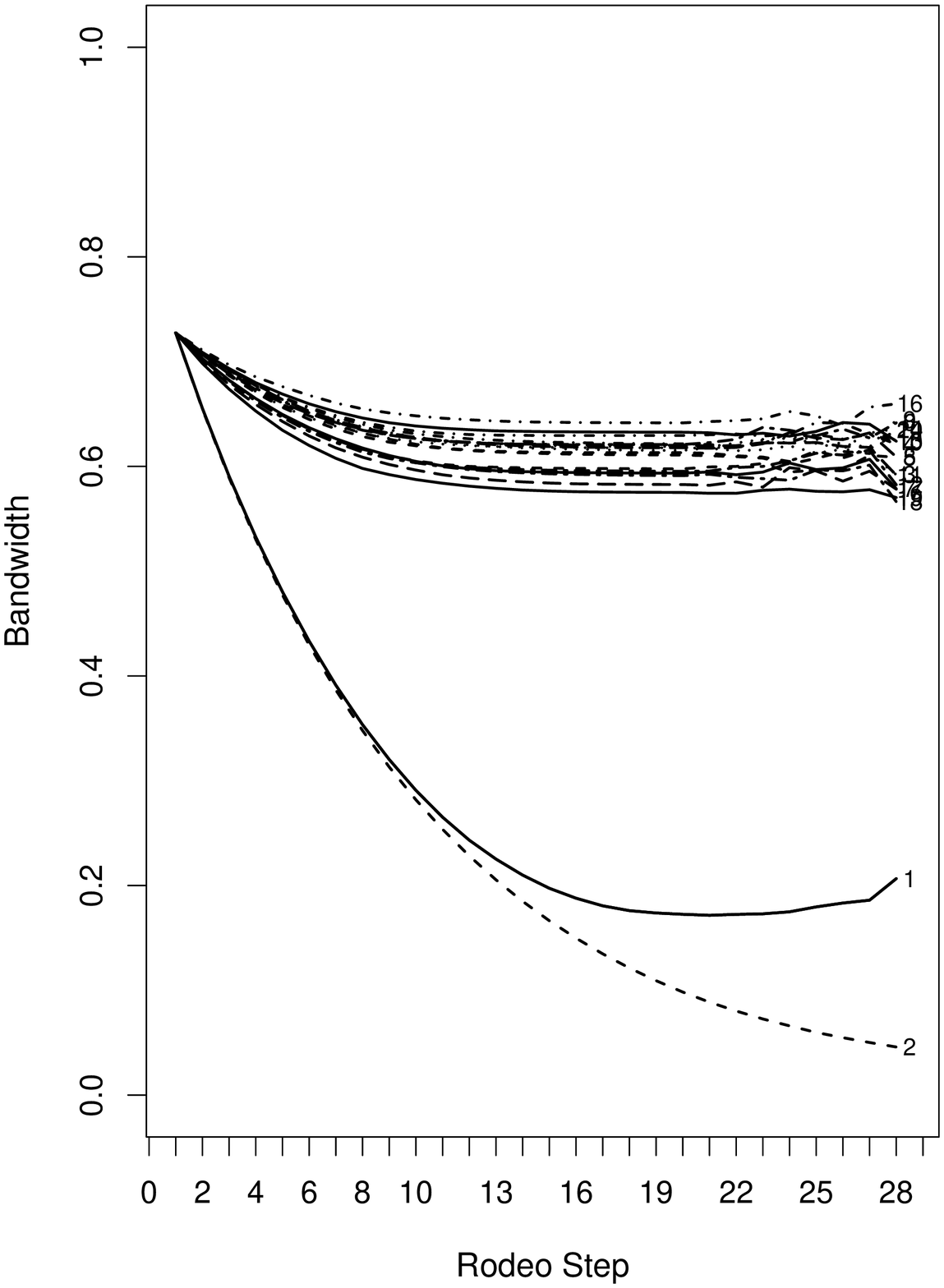} &
\includegraphics[width=2.0in,angle=0]{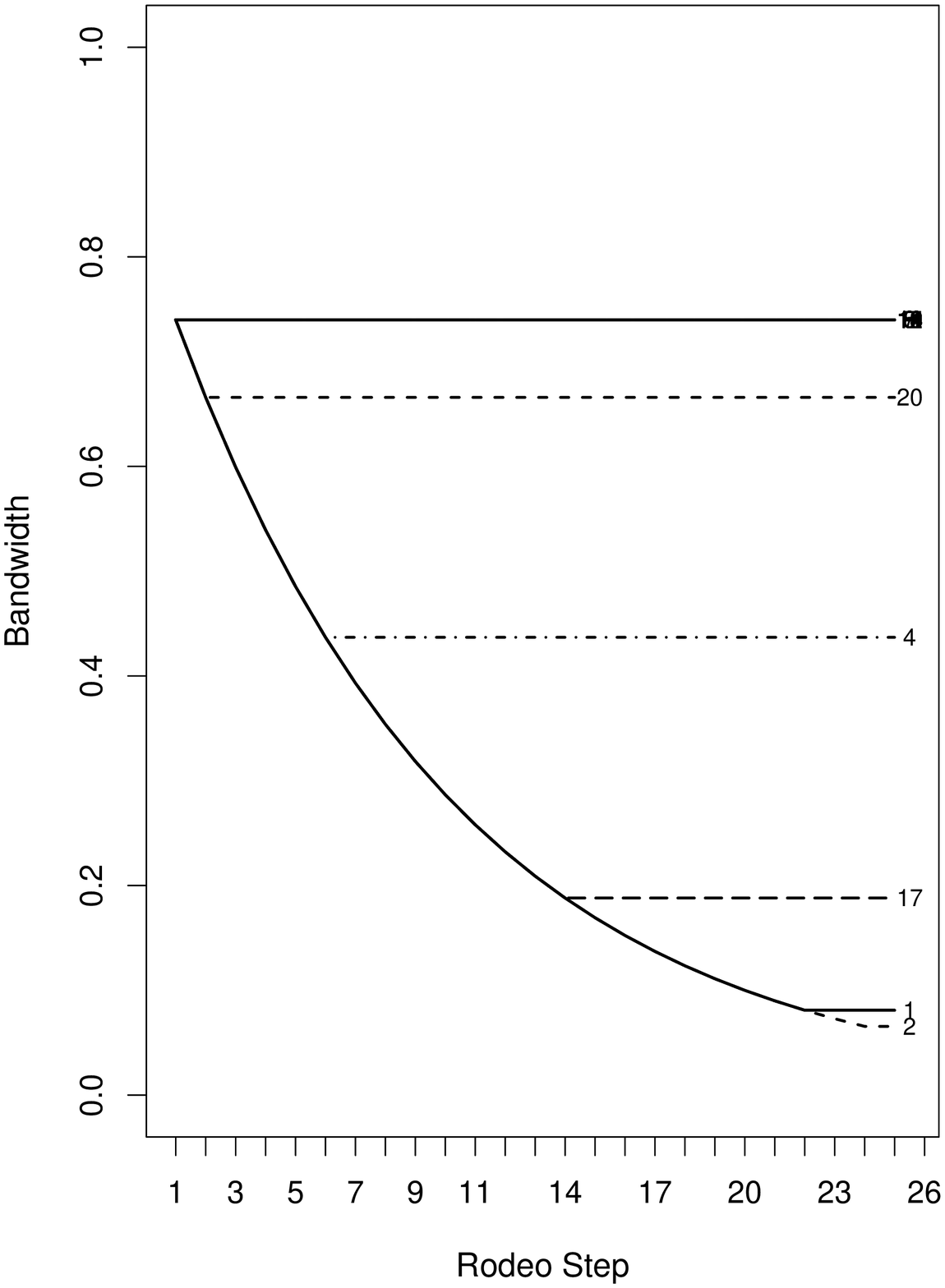} &
\includegraphics[width=2.0in,angle=0]{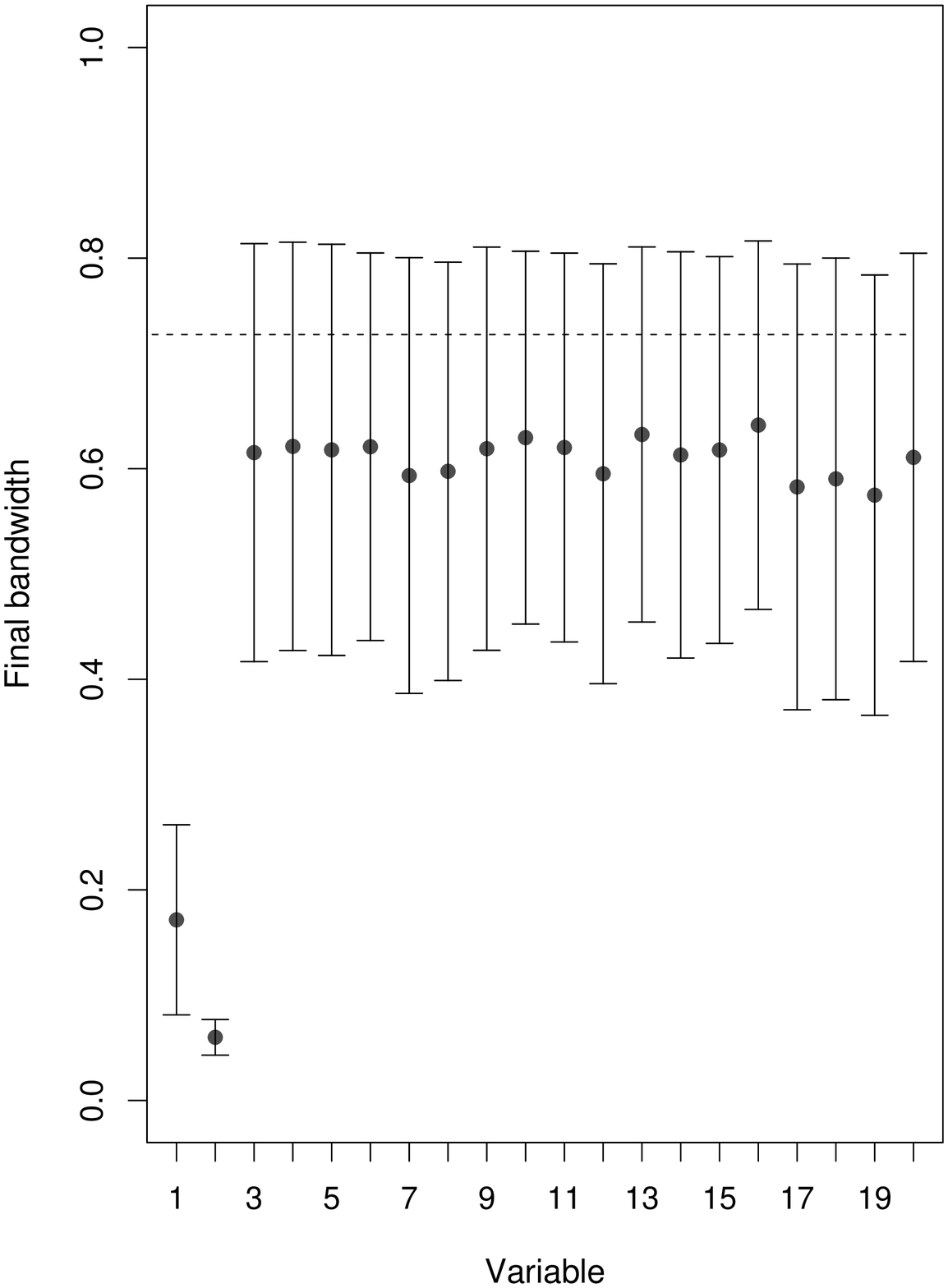} \\
\end{tabular}
\end{center}
\caption{Rodeo run on synthetic data sets, showing average
bandwidths over 200 runs (left), final bandwidths with
standard errors (right), and bandwidths on
a single run of the algorithm (center).
In the top plots the regression function is $m(x) = 5
x_1^2 x_2^2 $ with $d=10$, $n=500$, and $\sigma = .5$
and in the lower plots the regression function is $m(x) =
2(x_1+1)^3 + 2\sin(10x_2)$, $d=20$, $n=750$, and $\sigma = 1$.
The figures show that
the bandwidths for the relevant variables $x_1$ and $x_2$ are shrunk,
while the bandwidths for the irrelevant variables remain large.}
\label{fig::tworel}
\end{figure*}

In the first example, we take $m(x) = 5 x_1^2 x_2^2 $ with $d=10$,
$\sigma = .5$ with $x_i \sim \text{Uniform}(0,1)$.  The algorithm is
applied to the local linear estimates around the test point $x_0 =
(\frac{1}{2},\ldots, \frac{1}{2})$, with $\beta=0.8$.  Figure
\ref{fig::tworel} shows the bandwidths averaged over 200 runs of
the rodeo, on data sets of size $n=750$.  
The second example in Figure~\ref{fig::tworel} shows the
algorithm applied to the function $m(x) = 2(x_1+1)^3 + 2\sin(10x_2)$, in
this case in $d=20$ dimensions with $\sigma = 1$.  

The plots demonstrate how the
bandwidths $h_1$ and $h_2$ of the relevant variables are shrunk, while
the bandwidths of the relevant variables tend to remain large.

\begin{figure*}
\begin{center}
\begin{tabular}{cc}
\includegraphics[width=2.3in,angle=0]{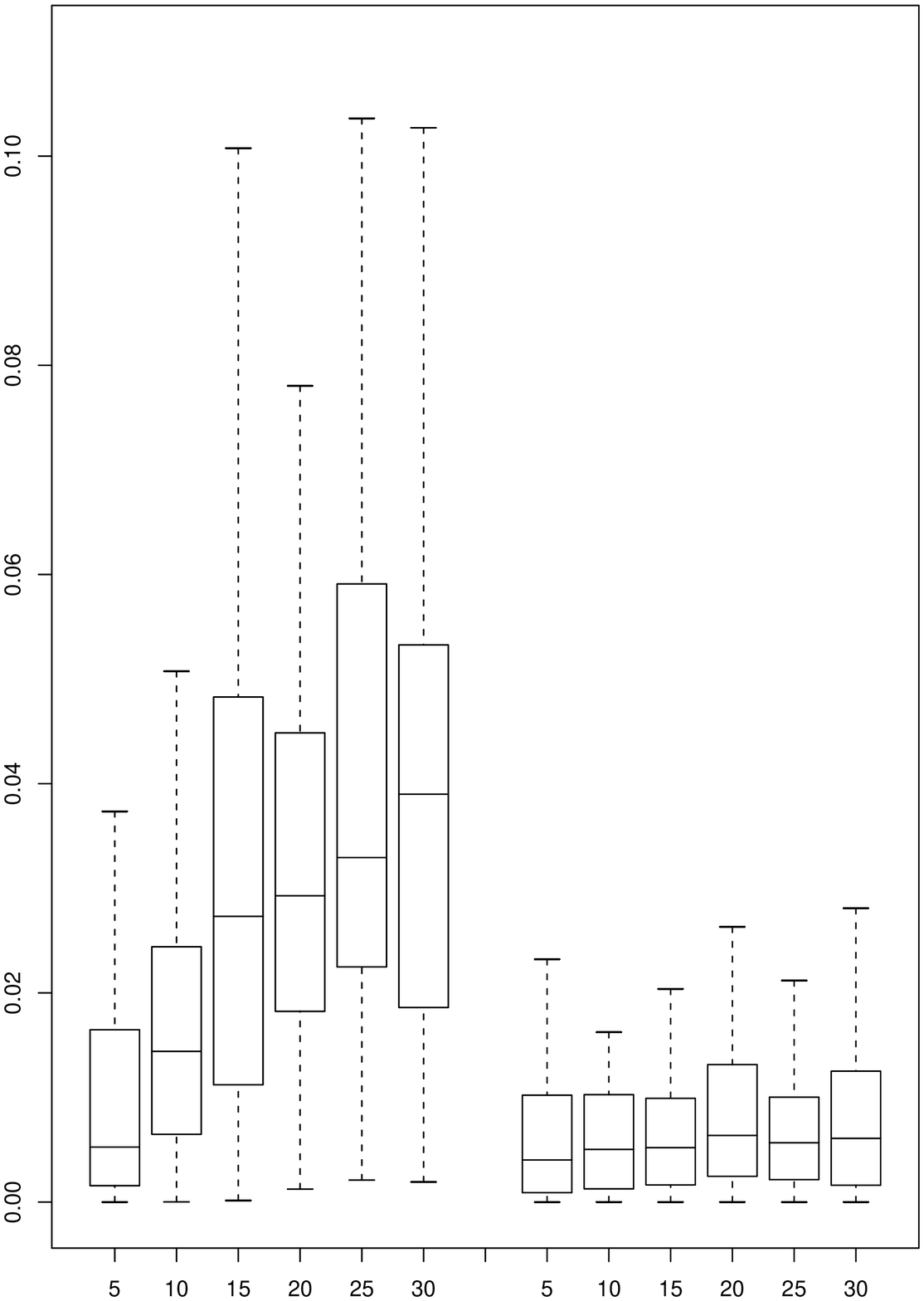} &
\includegraphics[width=2.3in,angle=0]{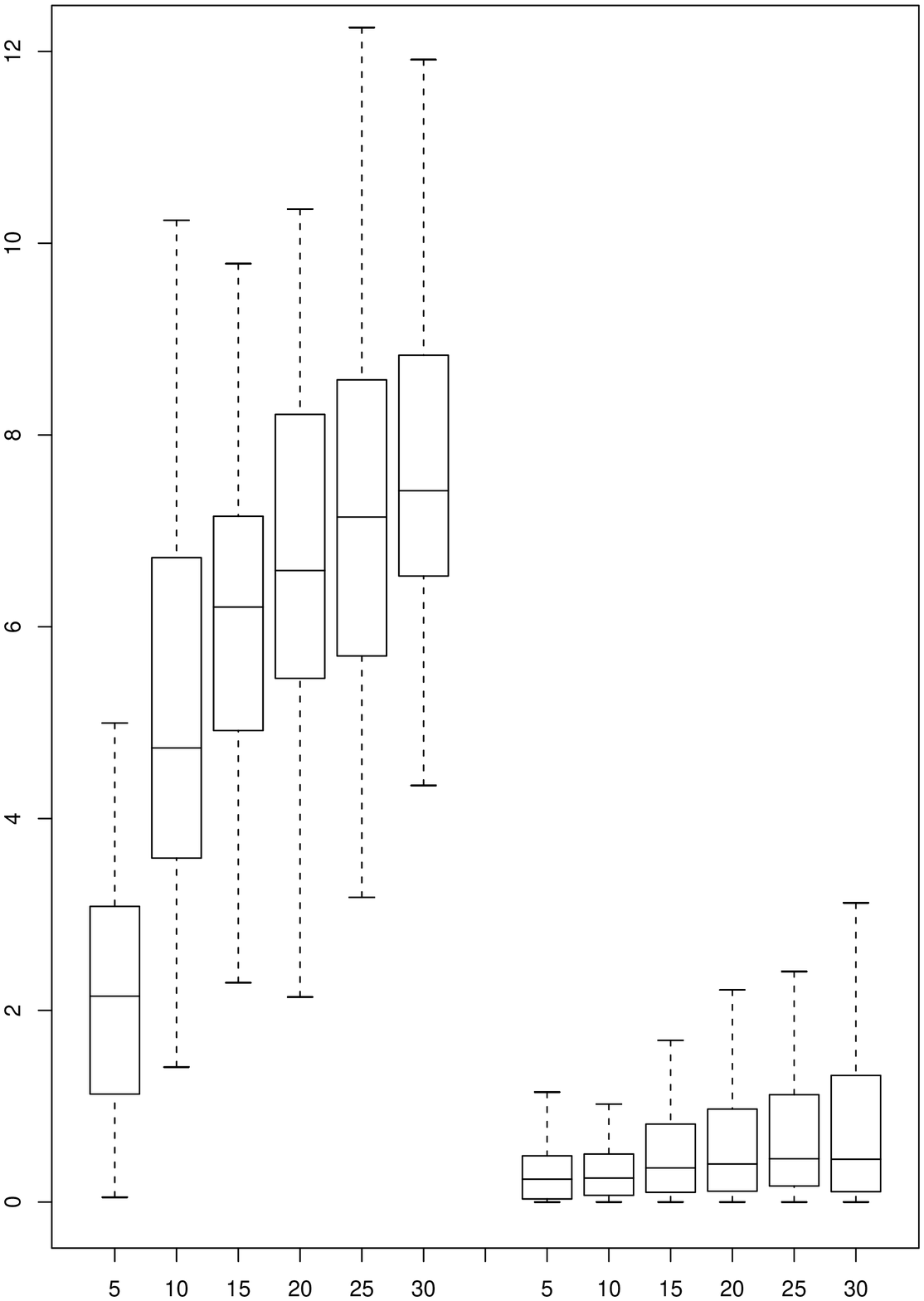} \\
\end{tabular}
\end{center}
\caption{Squared error of the estimator on the previous
examples, $m(x) = 5 x_1^2 x_2^2 $ (left) and $m(x) = 2(x_1+1)^3 +
2\sin(10x_2)$ right. For each plot, the left six boxplots show the
risk in different dimensions ($d=5,10,15,20,25,30$) when using a
single bandwidth, chosen by leave-one-out cross validation.  The right
six boxplots show the squared error on the same data with
bandwidths selected using the rodeo.}
\label{fig::risk}
\end{figure*}


\begin{figure}[t]
\begin{center}
\begin{tabular}{cc}
\includegraphics[width=2.3in,angle=0]{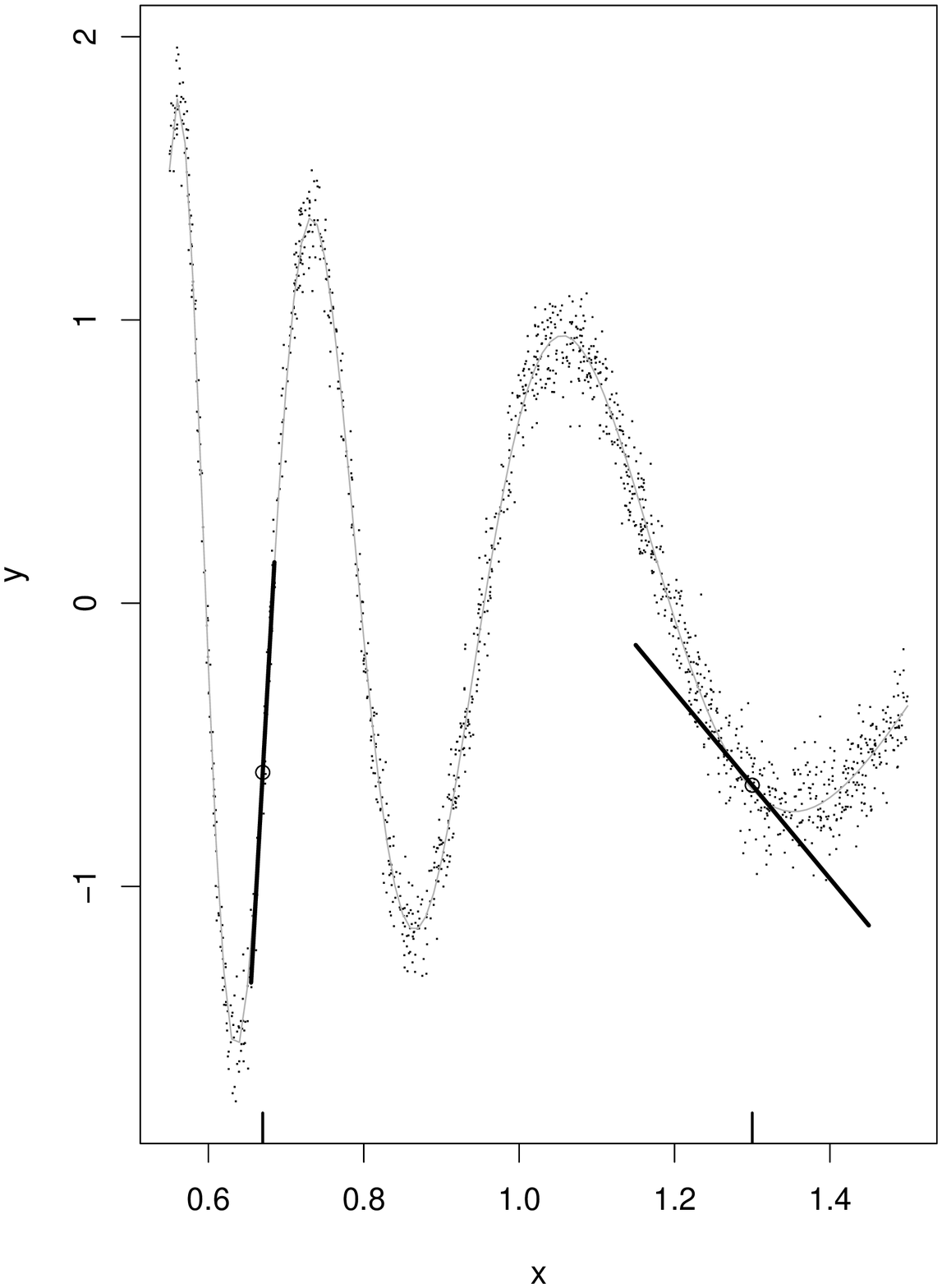} &
\includegraphics[width=2.3in,angle=0]{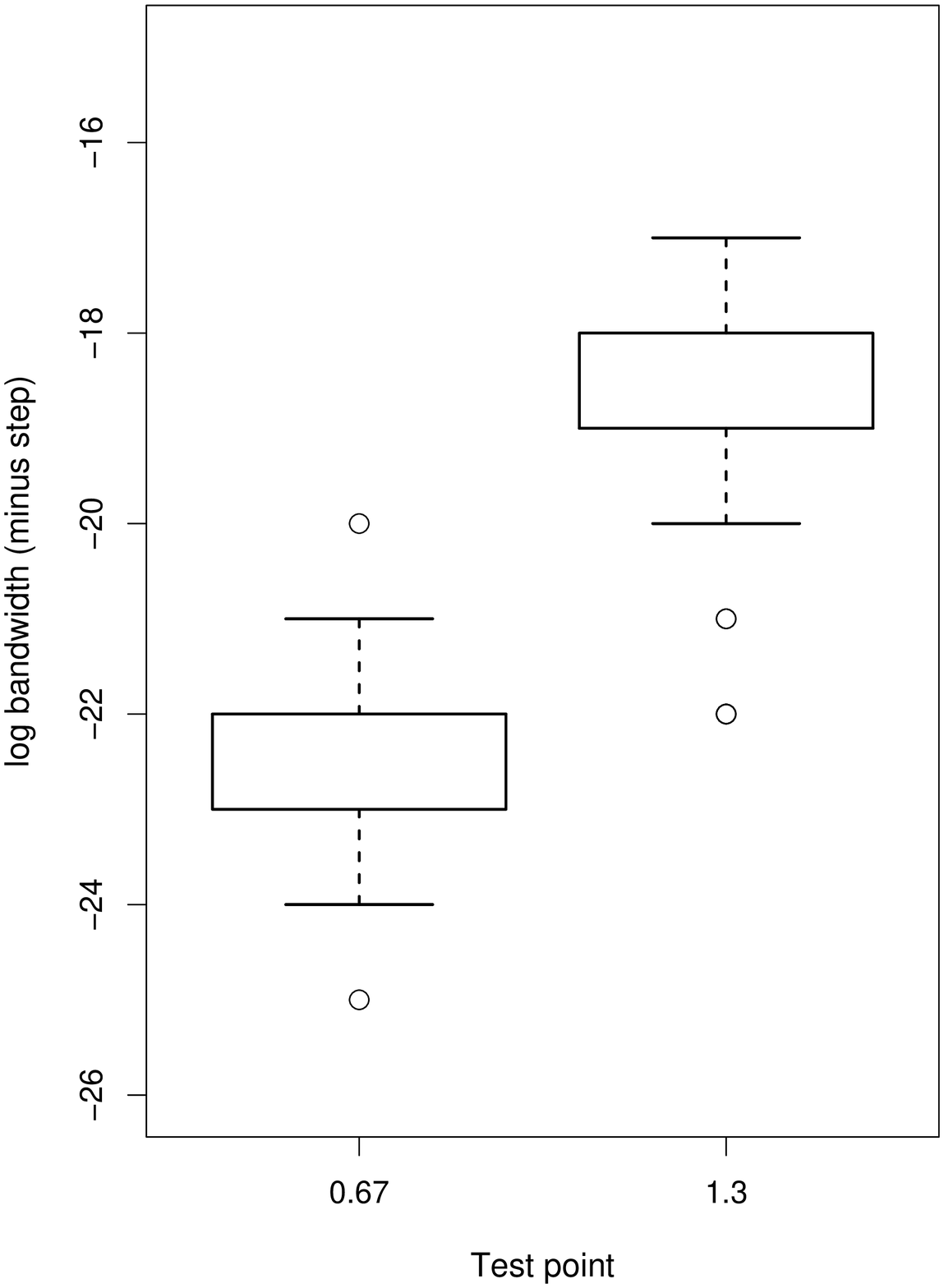}
\\[-.2in]
\end{tabular}
\end{center}
\caption{A one-dimensional example.  The regression
function is $m(x) = (1/x) \sin(15/x)$, and $n=1500$ data
points are sampled, $x\sim \text{Uniform}(0,1)+\frac{1}{2}$.  The left
plot shows the local linear fit at two test points; the 
right plot shows the final log bandwidth, $\log_{1/\beta} h_\star$, 
(equivalently, minus the number of steps)
of the rodeo over 50 randomly generated data sets.}
\label{fig::1dex}
\end{figure}

\subsection{A One-Dimensional Example}

The next figure illustrates the algorithm in one dimension.  The
underlying function in this case is $m(x) = (1/x) \sin(15/x)$, and
$n=1500$ data points are sampled as $x\sim
\text{Uniform}(0,1)+\frac{1}{2}$.  The algorithm is run at two test
points; the function is more rapidly varying near the test point
$x=0.67$ than near the test point $x=1.3$, and the rodeo
appropriately selects a smaller bandwidth at $x=0.67$.  The right
plot of Figure~\ref{fig::1dex} displays boxplots for 
logarithm of the final bandwidth, in the base $1/\beta$
(equivalently, minus the number of steps in the algorithm),
averaged over 50 randomly generated data sets.  

The figure illustrates how smaller bandwidths are selected where the
function is more rapidly varying.  However, we do not claim that the
method is adaptive over large classes of function spaces.  As
discussed earlier, the technique is intentionally a greedy algorithm;
adapting to unknown smoothness may require a more refined search over
bandwidths that does not scale to large dimensions, and is 
out of the scope of the current paper.

\subsection{Estimating $\sigma$}\label{sec::sigma}

The algorithm requires that we insert
an estimate $\hat\sigma$ of $\sigma$ in (\ref{eq::sj}).  An estimator
for $\sigma$ can be obtained by generalizing a method of
\cite{Rice:1984}.  For $i<\ell$, let 
\begin{equation}
d_{i\ell} = \|X_{i} - X_{\ell}\|\,.
\end{equation}  
Fix an integer $J$ and let ${\cal E}$ denote the set of
pairs $(i,\ell)$ corresponding the $J$ smallest values of $d_{i\ell}$.
Now define
\begin{equation}
\hat{\sigma}^2 = \frac{1}{2J}\sum_{i,\ell\in{\cal E}} (Y_{i} - Y_{\ell})^2.
\end{equation}
Then,
\begin{equation}
\E(\hat{\sigma}^2) \;=\; \sigma^2 + {\rm bias}
\end{equation}
where
\begin{equation}
{\rm bias} \;\leq\;  D  \, \sup_x \sum_{j=1}^r \left| \frac{\partial m_j(x)}{\partial x_j} \right| 
\end{equation}
with $D$ given by
\begin{equation}
D = \max_{i,\ell\in {\cal E}} \|X_{i} - X_{\ell}\|.
\end{equation}
There is a bias-variance tradeoff:
large $J$ makes $\hat{\sigma}^2$ positively biased,
and small $J$ makes $\hat{\sigma}^2$ highly variable.
Note, however, that the bias is mitigated by sparsity (small $r$).

A more robust estimate may result from taking
\begin{equation}
\label{eq::sigmahat}
\hat{\sigma}^2 = 
\frac{\sqrt{\pi}}{2} \mathop{\text{median}} \left\{\left|Y_{i} - Y_{\ell}\right|\right\}_{i,\ell\in{\cal E}} 
\end{equation}
where the constant comes from observing that if $X_i$ is close to
$X_{\ell}$, then 
\begin{equation}
|Y_i-Y_\ell| \sim |N(0,2\sigma^2)| = \sqrt{2}\sigma |Z|,
\end{equation}
where $Z$ is a standard normal with $\E|Z| = \sqrt{2/\pi}$.

\begin{figure*}
\begin{center}
\begin{tabular}{ccc}
& \\[-30pt]
&& \\[-30pt]
\includegraphics[width=2.0in,angle=0]{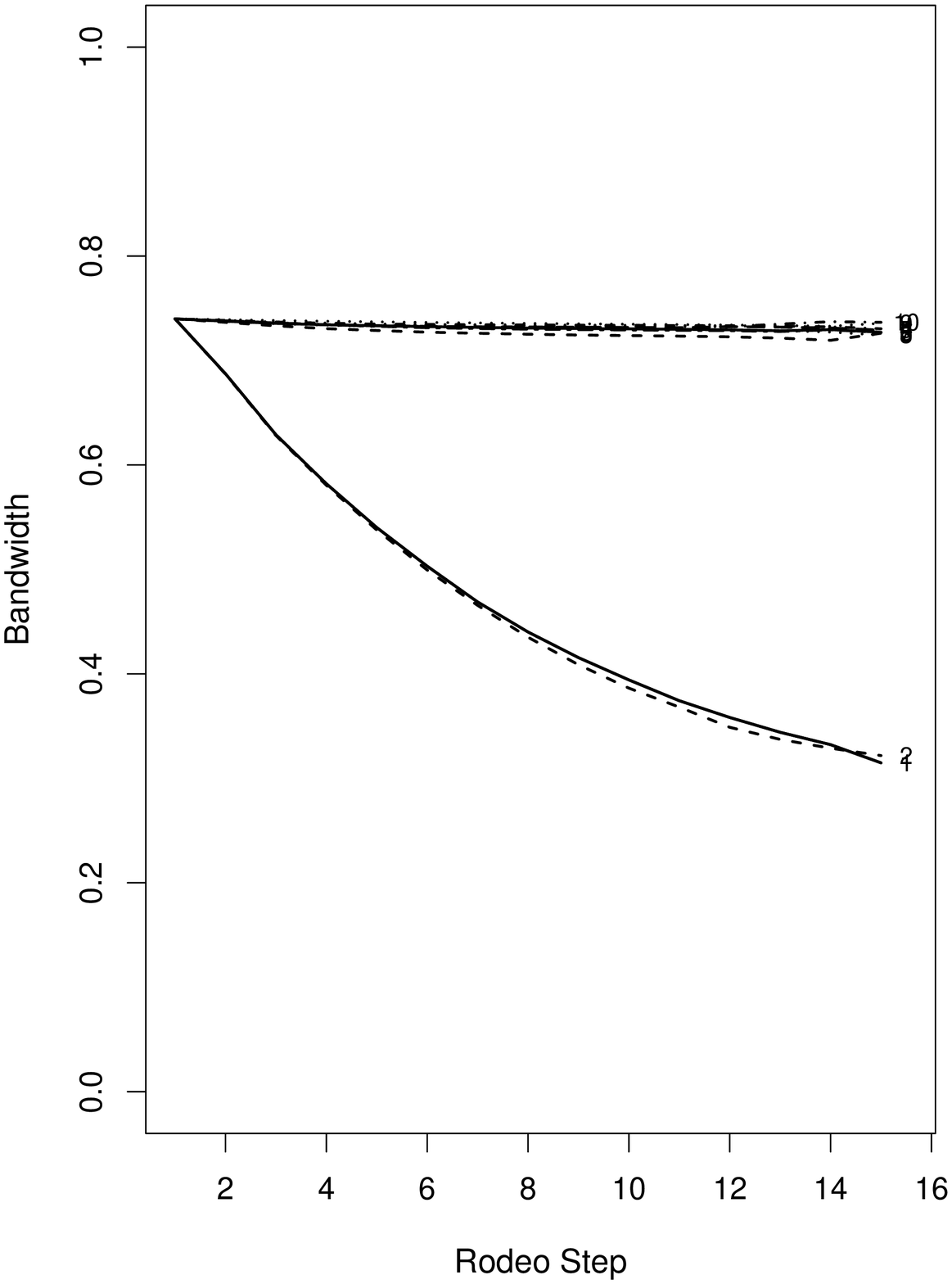} &
\includegraphics[width=2.0in,angle=0]{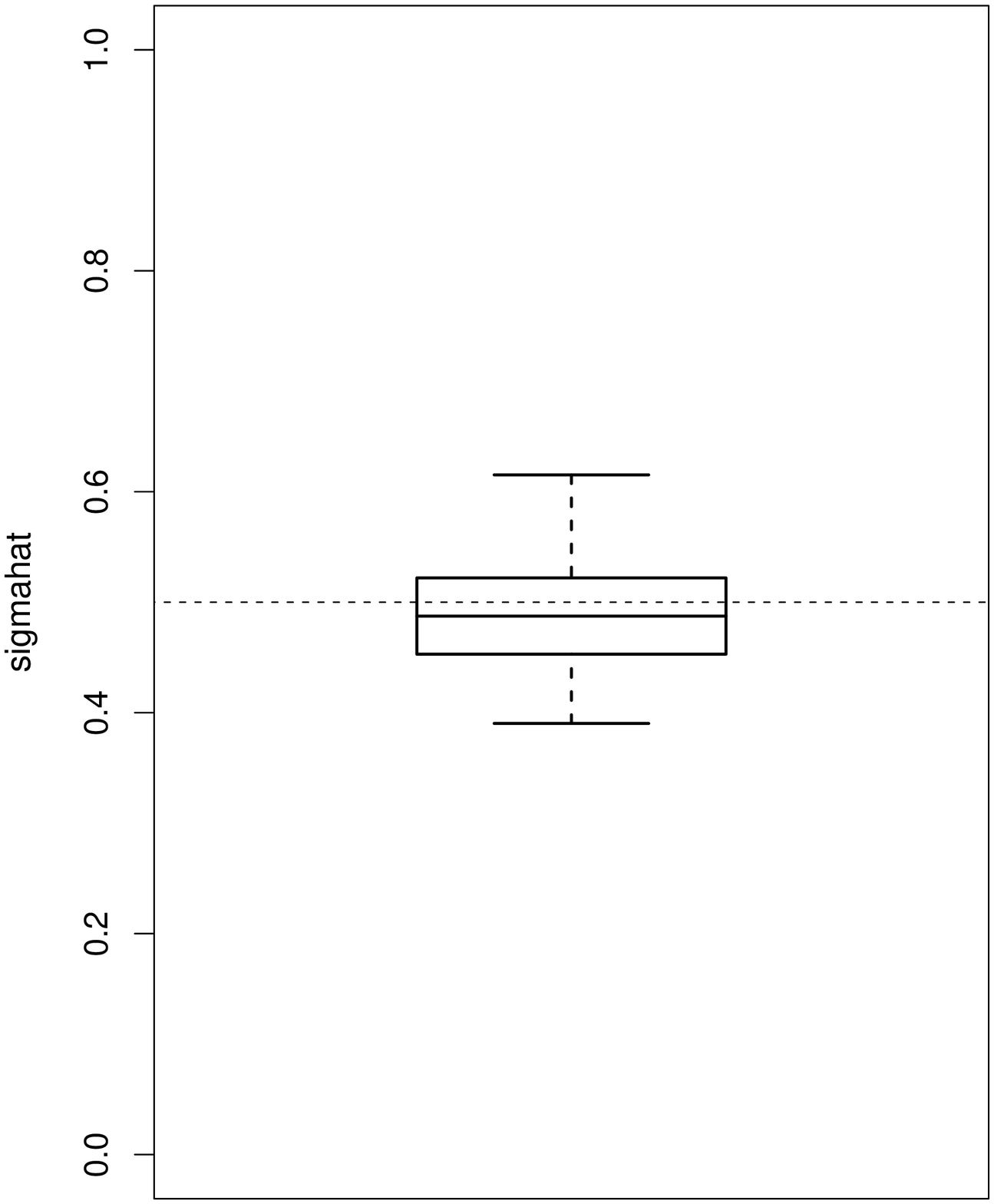} &
\includegraphics[width=2.0in,angle=0]{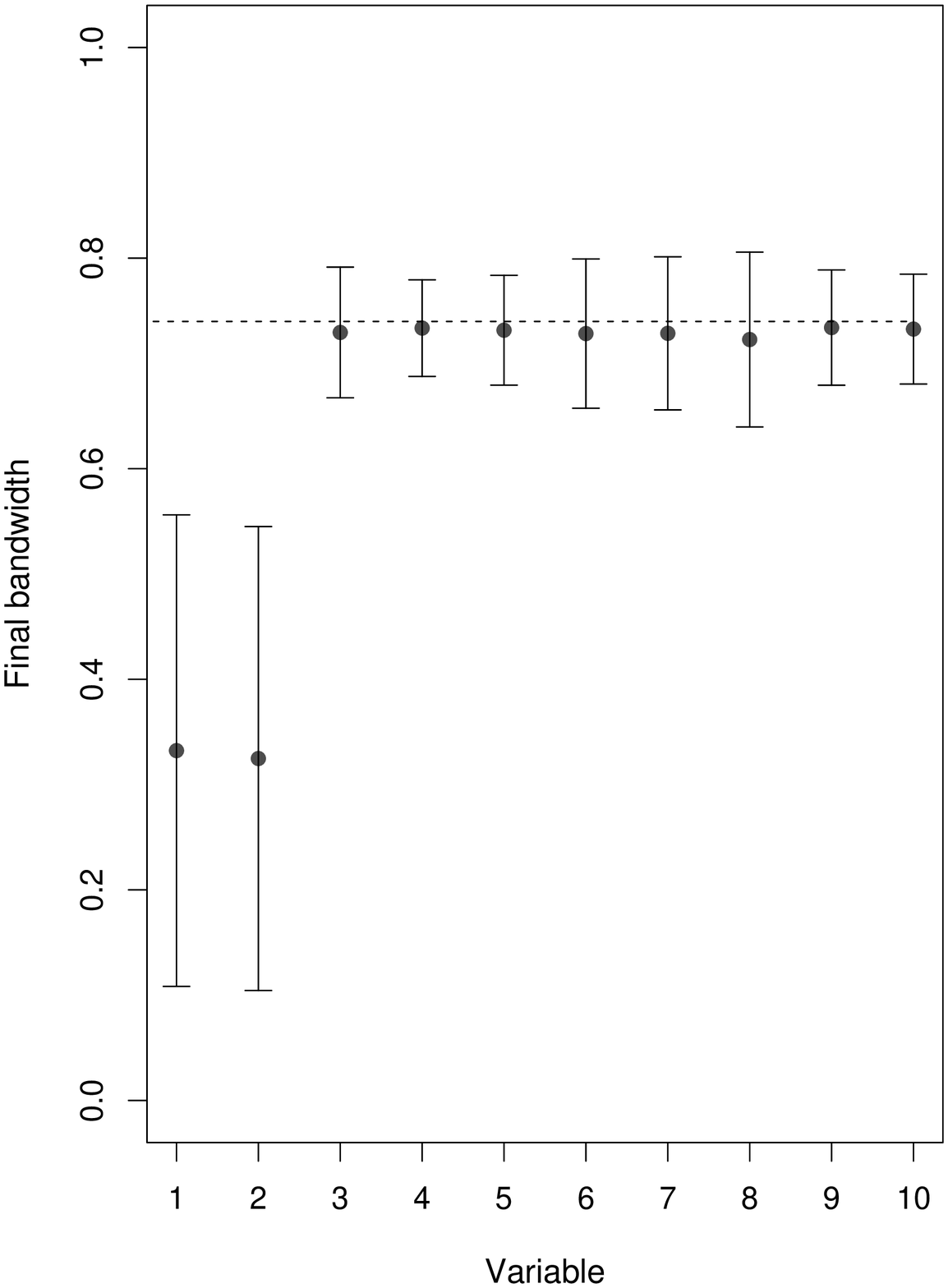} \\
&& \\[-30pt]
\includegraphics[width=2.0in,angle=0]{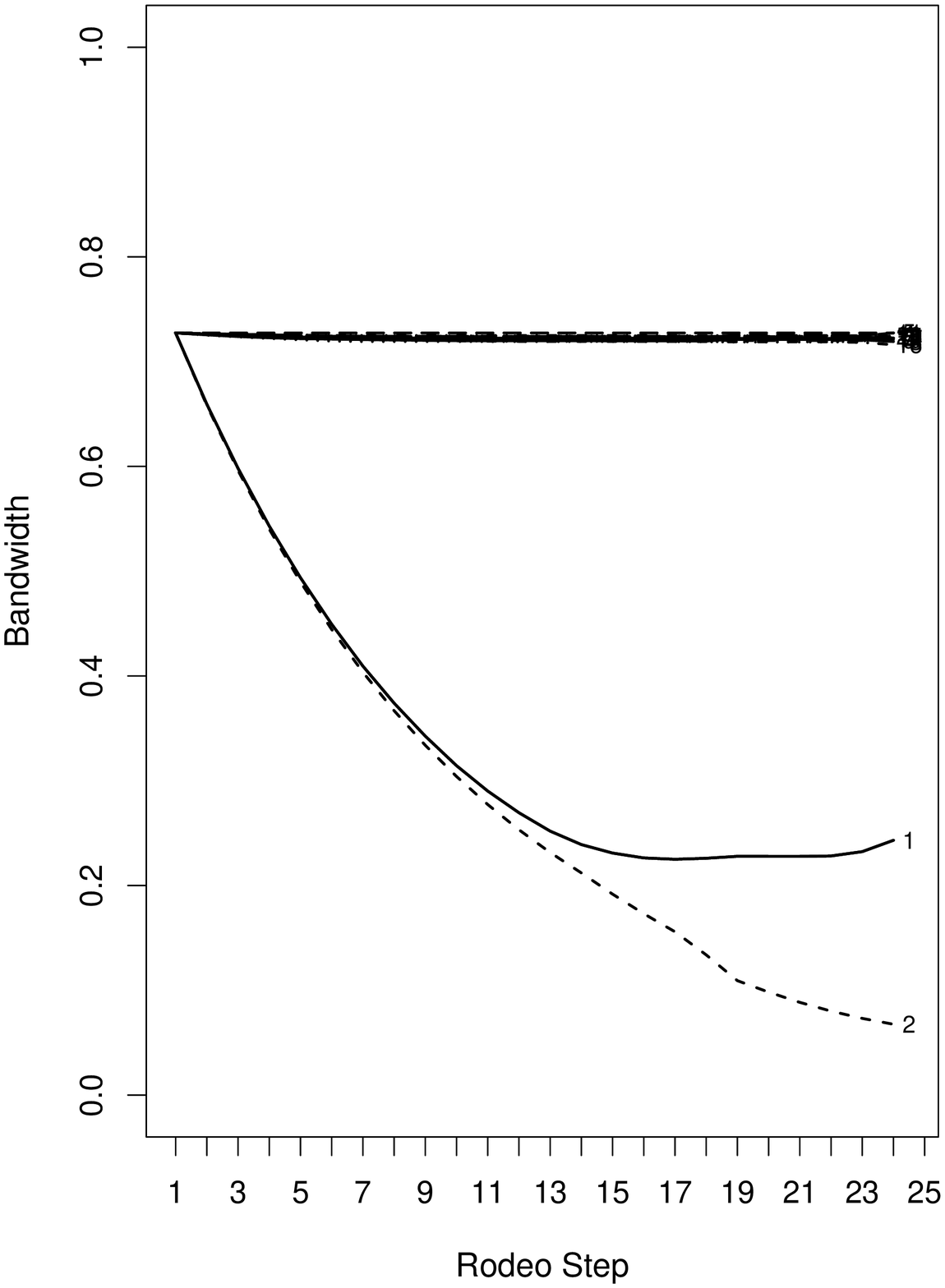} &
\includegraphics[width=2.0in,angle=0]{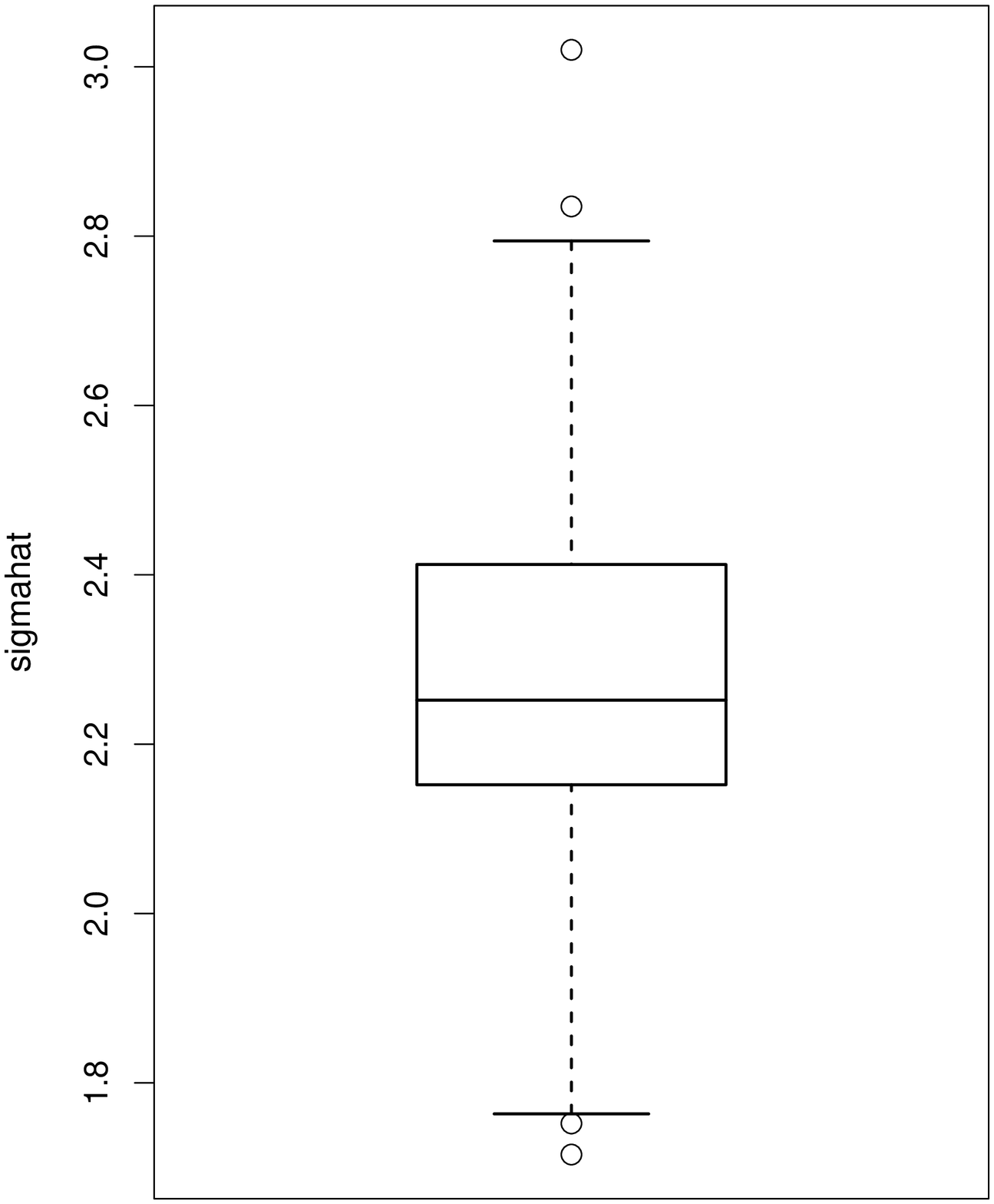} &
\includegraphics[width=2.0in,angle=0]{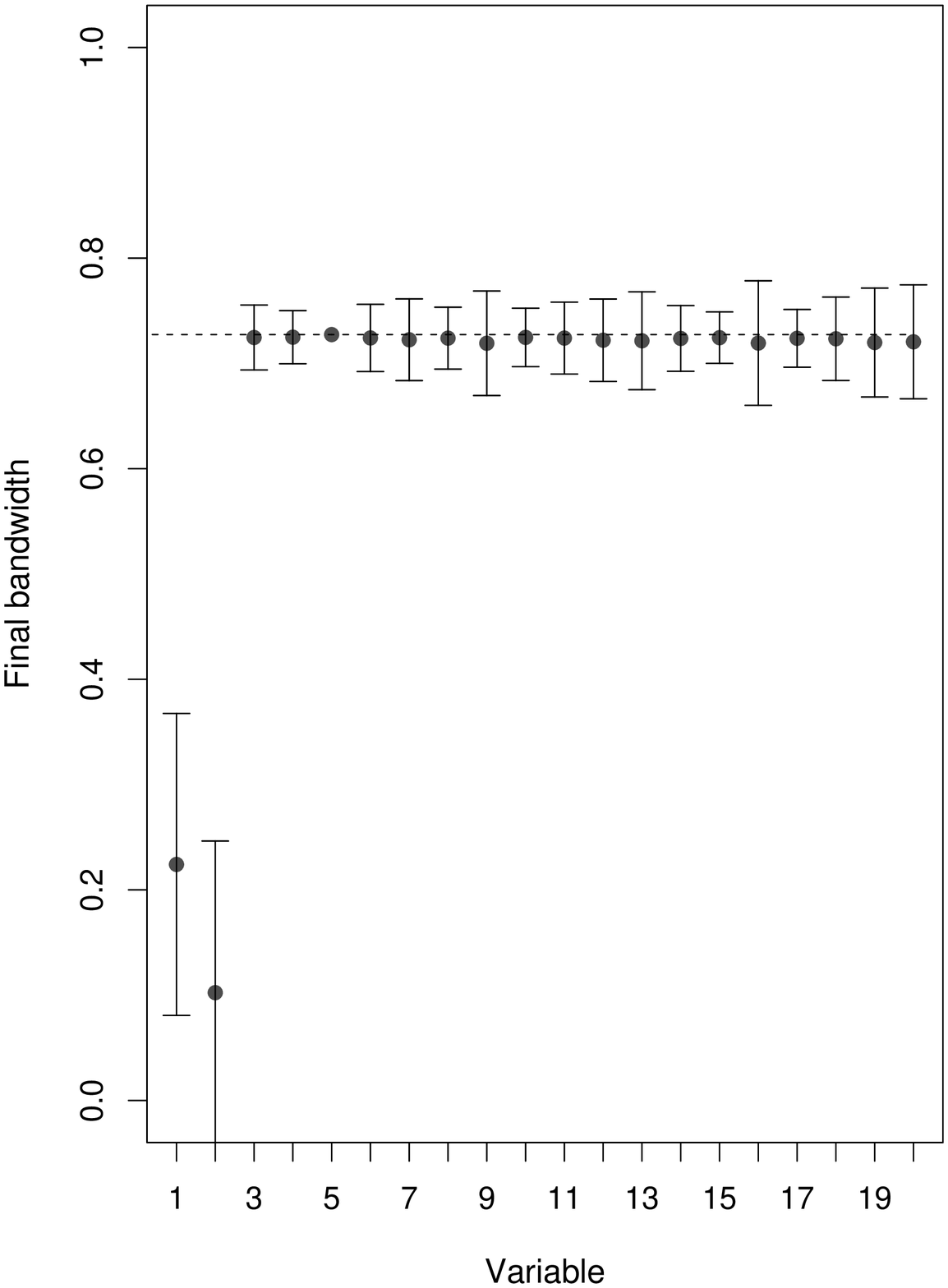} \\
\end{tabular}
\end{center}
\caption{Rodeo run on the examples of Section~\ref{sec::tworel},
but now estimating the noise using the estimate $\hat\sigma$ discussed in
Section~\ref{sec::sigma}. Top: $\sigma=.5$, $d=10$; bottom:
  $\sigma=1$, $d=20$.  In higher dimensions the noise is
over-estimated (center plots), which 
results in the irrelevant variables being more aggressively
  eliminated; compare Figure~\ref{fig::tworel}.}
\label{fig::sigmahatex}
\end{figure*}

Now we redo the earlier examples, taking $\sigma$ as unknown.
Figure~\ref{fig::sigmahatex} shows the result of running the
algorithm on the examples of Section~\ref{sec::tworel},
however now estimating the noise using estimate~\eqref{eq::sigmahat}.
For the higher dimensional example, with $d=20$, the noise
variance is over-estimated, with the primary result 
that the irrelevant variables are more aggressively thresholded
out; compare Figure~\ref{fig::sigmahatex} to Figure~\ref{fig::tworel}.

Although we do not pursue it in this paper,
there is also the possibility of allowing
$\sigma(x)$ to be a function of $x$ and estimating it locally.

\subsection{Computational Cost}

When based on a local linear estimator, each step of the rodeo
algorithm has the same computational cost as constructing a single
local linear fit.  This is dominated by the cost of constructing the
matrix inverse $(X^T W X)^{-1}$ in equation \eqref{eq:inverse}.  Since
the derivative needs to be computed for every variable, the algorithm
thus scales as $O(d^4)$ in the dimension $d$.  Implemented in R, the
20 dimensional example in Figure~\ref{fig::tworel} takes 4 hours, 4
minutes and 40 seconds for 200 runs, or 73.4 seconds per run, when
executed on a 1.5~GHz PowerPC Macintosh laptop.  Although we focus on
local linear regression, it should be noted that very similar results
are obtained with kernel regression, which requires no matrix
inversion.  Using kernel regression, the same example requires
12~minutes and 33~seconds, or 3.7 seconds per run.

\section{Properties of the Rodeo}
\label{sec::properties}

We now give some results on the statistical properties of the hard
thresholding version of the rodeo estimator.  Formally, we use a
triangular array approach so that $m(x)$, $f(x)$, $d$ and $r$ can all
change as $n$ changes, although we often suppress the dependence
on $n$ for notational clarity.  We assume throughout that $m$ has continuous
third order derivatives in a neighborhood of $x$.  For convenience of
notation we assume that the covariates are numbered such that the
relevant variables $x_j$ correspond to $1\leq j\leq r$ and the
irrelevant variables $x_j$ correspond to $r+1\leq j\leq d$.

A key aspect of our analysis is that we allow the dimension $d$ to
increase with sample size $n$, and show that the algorithm achieves
near optimal minimax rates of convergence if $d = O(\log n /\log\log
n)$.  This hinges on a careful analysis of the asymptotic bias and
variance of the estimated derivative $Z_j$, taking the increasing
dimension into account.  We conjecture that, without further
assumptions, $d$ cannot increase at a significantly faster rate,
while obtaining near optimal rates of convergence.

The results are stated below, with the
complete proofs given in Section~\ref{sec::proofs}.

We write
$Y_n=\tilde O_P(a_n)$ to mean that $Y_n=O_P(b_n
a_n)$ where $b_n$ is logarithmic in $n$.  As noted earlier,
we write $a_n = \Omega(b_n)$ if 
$\liminf_{n} \left|\frac{a_n}{b_n}\right| > 0$; similarly
$a_n = \tilde\Omega(b_n)$ if $a_n = \Omega(b_n c_n)$ where
$c_n$ is logarithmic in $n$.

Let
$
\H =
\begin{pmatrix}
\H_R & 0\\
0 & 0
\end{pmatrix}
$ 
denote the Hessian of $m(x)$,
let $h_j^{(t)}$ denote the $j^{\rm th}$ bandwidth at step $t$ 
of the algorithm and
denote the bandwidth matrix by $H^{(t)} = \text{diag}((h_1^{(t)})^2,\ldots,(h_d^{(t)})^2)$.
Similarly, let $H_R^{(t)} = \text{diag}((h_1^{(t)})^2,\ldots,(h_r^{(t)})^2)$.
Define
\begin{equation}
\mu_j^{(t)}=\frac{\partial}{\partial h_j} \E\left[\hat m_{H^{(t)}}(x)
  - m(x) \,|\, X_1,\ldots, X_n\right].
\end{equation}
Our first result is an asymptotic expansion of the bias
of the estimator.

\begin{theorem}\label{thm:important}
Suppose that $x$ is interior to the support of $f$ and
that $K$ is a product kernel.
Assume that
\begin{equation}
d  = O(\log n / \log\log n),\qquad r=O(1),
\end{equation}
\begin{equation}\label{eq::cond-on-h}
h_0 = \frac{c_0}{\log\log n}
\end{equation}
for some $c_0 >0$
and
\begin{equation}
\label{eq::alpha}
\beta_n = 
n^{-\alpha/\log^3 n}
\end{equation}
for some $0 < \alpha < 1$.
Define
\begin{equation}
L_j^{(t)} =
\left\{
\begin{array}{ll}
\displaystyle {\nu_2 m_{jj}(x)} h_j^{(t)} & j \leq r\\
\noalign{\vskip 4pt}
-\tr\left(H_R^{(t)} \H_R^{(t)}\right)\, \nu_2^{2} \left(\nabla_j \log f(x)\right)^2 h_j^{(t)} & j > r.
\end{array}
\right.
\end{equation}
where $\nu_2$ is defined in equation \eqref{eq:nudef}.
Suppose that
\begin{equation}\label{eq::d2}
\max_{a,b}\sup_x \left| \frac{\partial^2 f(x)}{\partial x_a\partial x_b}\right| =
O(1/d^2)
\end{equation}
and
\begin{equation}\label{eq::d3}
\max_{a,b,c}\sup_x \left| \frac{\partial^3 f(x)}
     {\partial x_a\partial x_b\partial x_c}\right|=
O(1/d^3).
\end{equation}
Then, for $T_n \leq c_1\log n$,
\begin{equation}
\max_{\stackrel{1\leq s\leq T_n}{ 1\leq j \leq d}}|\mu_j^{(t)} - L_j^{(t)}| =
O_P\left((h_j^{(t)})^2/d\right).
\end{equation}
\end{theorem}

The next result characterizes the variance
$v_j^{(t)} = {\rm Var}(Z_j(t) \,|\, X_1, \ldots, X_n)$.

\begin{theorem}
\label{thm:var}
Assume the conditions of Theorem~\ref{thm:important}.
Let
\begin{equation}
s_j^2(t)  = \frac{C}{n (h_j^{(t)})^2}\prod_{k=1}^d \frac{1}{h_k^{(t)}}
\end{equation}
where $h_j^{(t)}$ is the $j^{\rm th}$ bandwidth at step $t$,
with $C = {\sigma^2 \int K^2(u)\, du}/{f(x)}$.
Then,
\begin{equation}
\mathbb{P}\Biggl( 
\max_{\stackrel{1\leq s\leq T_n}{ 1\leq j \leq d}} 
 \left| \frac{v_j^{(t)}}{s_j^2(t)} -1 \right|  > \epsilon \Biggr) \to 0,
\end{equation}
for all $\epsilon >0$.
\end{theorem}

\vskip10pt
Our main theoretical result characterizes the asymptotic running time, 
selected bandwidths, and risk of the algorithm.  In order
to get a practical algorithm, we need to make assumptions on the
functions $m$ and $f$.

\begin{enumerate}
\item[(A1)]
For each $j > r$, 
\begin{equation}
\nabla_j \log f(x) = O\left( \frac{1}{n^{1/4}}\right).
\end{equation}
\item[(A2)] 
\begin{equation}
\label{eq::A}
\liminf_n \min_{1\leq j\leq r} |m_{jj}(x)| > 0 \,.
\end{equation} 
\end{enumerate}

\noindent
{\it Explanation of the Assumptions.}
To give the intuition behind these assumptions,
recall from Theorem \ref{thm:important} that
\begin{equation}\label{eq::AndB}
\mu_j = 
\left\{
\begin{array}{ll}
A_j h_j + o_P(h_j) & j \leq r\\
B_j h_j + o_P(h_j) & j > r
\end{array}
\right.
\end{equation}
where
\begin{equation}
A_j = {\nu_2 m_{jj}(x)},\ \ \ \ 
B_j = -\tr (H_R {\cal H}_R)\nu_2^2 (\nabla_j \log f(x))^2 .
\end{equation}
Moreover, 
$\mu_j \approx 0$ when the sampling density $f$ is uniform or the data are on a regular grid.
Consider assumption (A1).
If $f$ is uniform then
this assumption is automatically satisfied
since then ${\mu}_j(s)\approx 0$ for $j>r$.
More generally,
$\mu_j$ is approximately proportional to
$(\nabla_j \log f(x))^2$ for $j>r$ which implies that
$|\mu_j| \approx  0$ for irrelevant variables
if $f$ is sufficiently smooth in the variable $x_j$.
Hence, assumption (A1) can be interpreted as requiring that $f$ is
sufficiently smooth in the irrelevant dimensions.  
If $d$ is constant and $f$ is smooth,
then this assumption is \emph{implied} by Theorem ~\ref{thm:important}.
However, we allow for the possibility
that $d$ grows with $n$.

Now consider assumption (A2).
Equation (\ref{eq::AndB}) ensures that
$\mu_j$ is proportional to $h_j |m_{jj}(x)|$
for small $h_j$.  Since we take the initial bandwidth $\LL$ 
to be decreasingly slowly with $n$, 
(A2) implies that $|\mu_j(h)| \geq c h_j |m_{jj}(x)|$ for
some constant $c > 0$, for sufficiently large $n$.

\vskip20pt

\begin{theorem}
Assume the conditions of Theorem~\ref{thm:important}
and suppose that assumptions (A1) and (A2) hold. In addition, suppose
that 
$\Amin = \min_{j\leq r} |m_{jj}(x)| =
\tilde\Omega(1)$ and $\Amax = \max_{j\leq r}
|m_{jj}(x)| = \tilde O(1)$.
Then the rodeo outputs bandwidths $h^\star$ that satisfy
\begin{equation}
\label{eq:bw1}
\P\left(h^\star_j \,=\, \LL \; \mbox{for all $j >
  r$}\right) \;\longrightarrow 1 
\end{equation}
and for every $\epsilon > 0$, 
\begin{equation}
\label{eq:bw2}
\P\left( n^{-1/(4+r)-\epsilon} \,\leq\,  h^\star_j \, \leq\,
         n^{-1/(4+r)+\epsilon} 
\;\mbox{for all $j \leq r$}\right) \longrightarrow 1\,.
\end{equation}
\label{thm}
\end{theorem}

\begin{corollary}
Under the conditions of Theorem~\ref{thm}, the risk
of the rodeo estimator satisfies
\begin{equation}
\R(h^\star) \;=\; \E\int (\widehat{m}_{h^\star}(x) -
m(x))^2 \,dx \;=\; {O}_P\left(n^{-4/(4+r)+\epsilon}\right)
\end{equation} 
for every $\epsilon > 0$.
\label{cor}
\end{corollary}

The corollary follows easily from the previous results
using the bias-variance decomposition.  Details
are given in Section~\ref{sec::proofs}.
In case $d=O(1)$, then
the result can be strengthened to
\begin{equation}
\R(h^\star) \;=\; \widetilde{O}_P\left(n^{-4/(4+r)}\right)
\end{equation}
so that the optimal rate of convergence is obtained
up to logarithmic factors, as if the relevant variables
were known in advance and isolated.

Assumption (A1) imposes a condition on the sampling
density of the irrelevant variables, requiring that
it is sufficiently smooth in the irrelevant dimensions.  
Using Theorem~\ref{thm:important}, we can modify the
threshold used by the algorithm in order to
remove this assumption.

\begin{theorem}
Consider a modified rodeo algorithm where the threshold is
changed to 
\begin{equation}
\lambda_j = \rho_n\beta^{3t} + s_j\sqrt{2\log(n c_n)}
\end{equation}
at step $t$, where $\rho_n >0$ satisfies
\begin{equation}\label{eq::rho}
\frac{\rho_n}{h_0^3} \longrightarrow \infty,\quad \frac{\rho_n}{h_0}
\longrightarrow 0\quad \mbox{as $n\longrightarrow\infty$}.
\end{equation}
Then the results of Theorem~\ref{thm} hold under
assumption (A2), without assuming~(A1).
\label{thm:noassumption}
\end{theorem}

\section{Extensions and Variations of the Rodeo}
\label{sec::mods}

The rodeo represents a general strategy for nonparametric estimation,
based on the idea of regularizing or testing the derivatives of an
estimator with respect to smoothing parameters.  There are many ways
in which this basic strategy can be realized.  In this section we
discuss several variants of the basic hard thresholding version of the
rodeo, including a soft thresholding version, a global rather than
local bandwidth selection procedure, the use of testing and
generalized cross validation, and connections to least angle
regression.  Further numerical examples are also given to illustrate
these ideas.

\subsection{Subtracting off a Linear Lasso}\label{sec::subtract}

Local linear regression is a nonparametric method that contains linear
regression as a special case when $h\to\infty$.  If the true function
is linear but only a subset of the variables are relevant, then the
rodeo will fail to separate the relevant and irrelevant variables
since relevance is defined in terms of departures from the limiting
parametric model.  Indeed, the results depend on the Hessian of $m$
which is zero in the linear case.  The rodeo may return a full linear
fit with all variables.  A simple modification fixes this problem.
First, do linear variable selection using, say, the lasso
\citep{Tibshirani:96}.  Then run the rodeo on the residuals from that
fit, but using all of the variables.  An example of this procedure is
given below in Section~\ref{sec::lars}.

\subsection{Other Estimators and Other Paths}


\begin{figure*}
{\it Rodeo: Soft thresholding version}
\vskip5pt
\hrule
\vskip5pt
\begin{enumerate}
\item {\it Select\/} parameter $0 < \beta < 1$ and initial bandwidth $\LL$,
satisfying $1 \leq \LL \leq \log^{\ell/d} n$, for a fixed constant
$\ell$.   Let $c_n$ be a sequence satisfying $c_n = O(1)$.
\item {\it Initialize\/} the bandwidths, and activate all covariates:
\begin{enumerate}
\item $h_j = \LL$, $j=1,2,\ldots, d$.
\item $\A = \{1,2,\ldots, d\}$
\item Initialize step, $t=1$.
\end{enumerate}
\item {\it While $\A$ is nonempty}
\begin{enumerate}
\item Set $dh_j(t)=0$, $j=1, \ldots, d$.
\item Do for each $j\in \A$:
\begin{enumerate}
\itemsep=5pt
\item[(1)] Compute the estimated derivative expectation $Z_j$ and $s_j$.
\item[(2)] Compute the threshold
$ \lambda_j = \displaystyle s_j \sqrt{2 \log (n c_n)}$.
\item[(3)] If $|Z_j| > \lambda_j$, 
   set $dh_j(t) = (1-\beta)\, h_j$  and $h_j \leftarrow \beta h_j$; \\
  otherwise remove $j$ from $\A$.
\item[(4)] Set $\hat{D}_j(t) = {\rm sign}(Z_j(h))\,( |Z_j(h)| - \lambda_j)_+$.
\end{enumerate}
\item Increment step, $t \leftarrow t+1$.
\end{enumerate}
\item {\it Output} bandwidths $h^\star=(h_1,\ldots,h_d)$ and estimator
\begin{equation}
\tilde{m}(x)  =  \hat{m}_{\LL}(x) - \sum_{s=1}^t  \bigl\langle
\hat{D}(s), dh(s)\bigr\rangle
\end{equation}
\end{enumerate}
\hrule
\vskip5pt
\caption{The soft thresholding version of the rodeo.}
\label{fig::soft-rodeo}
\vskip15pt
\end{figure*}

\begin{figure}[t]
\begin{center}
\begin{tabular}{ccc}
\includegraphics[width=2.3in,angle=0]{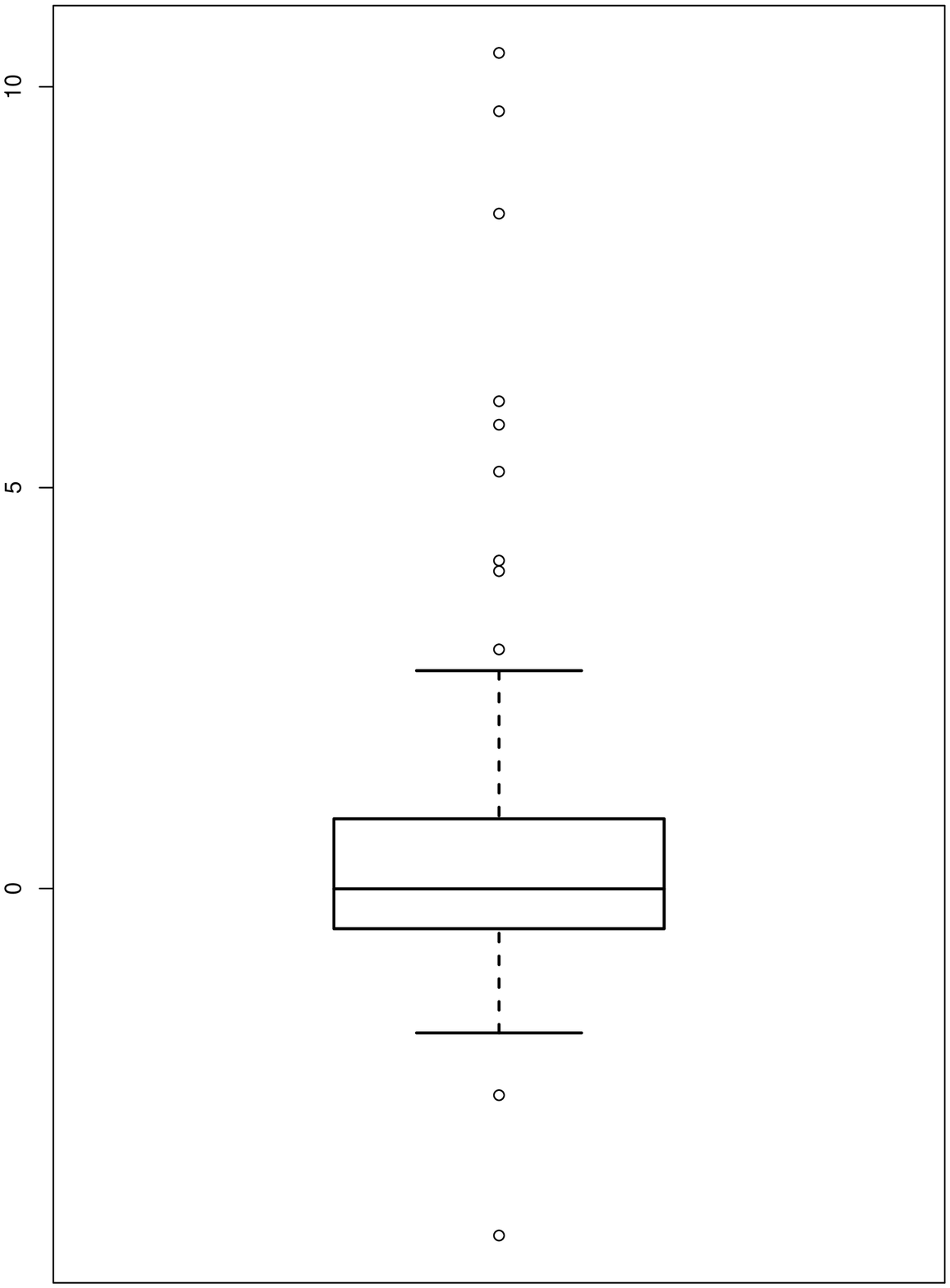}
&\quad&
\includegraphics[width=2.3in,angle=0]{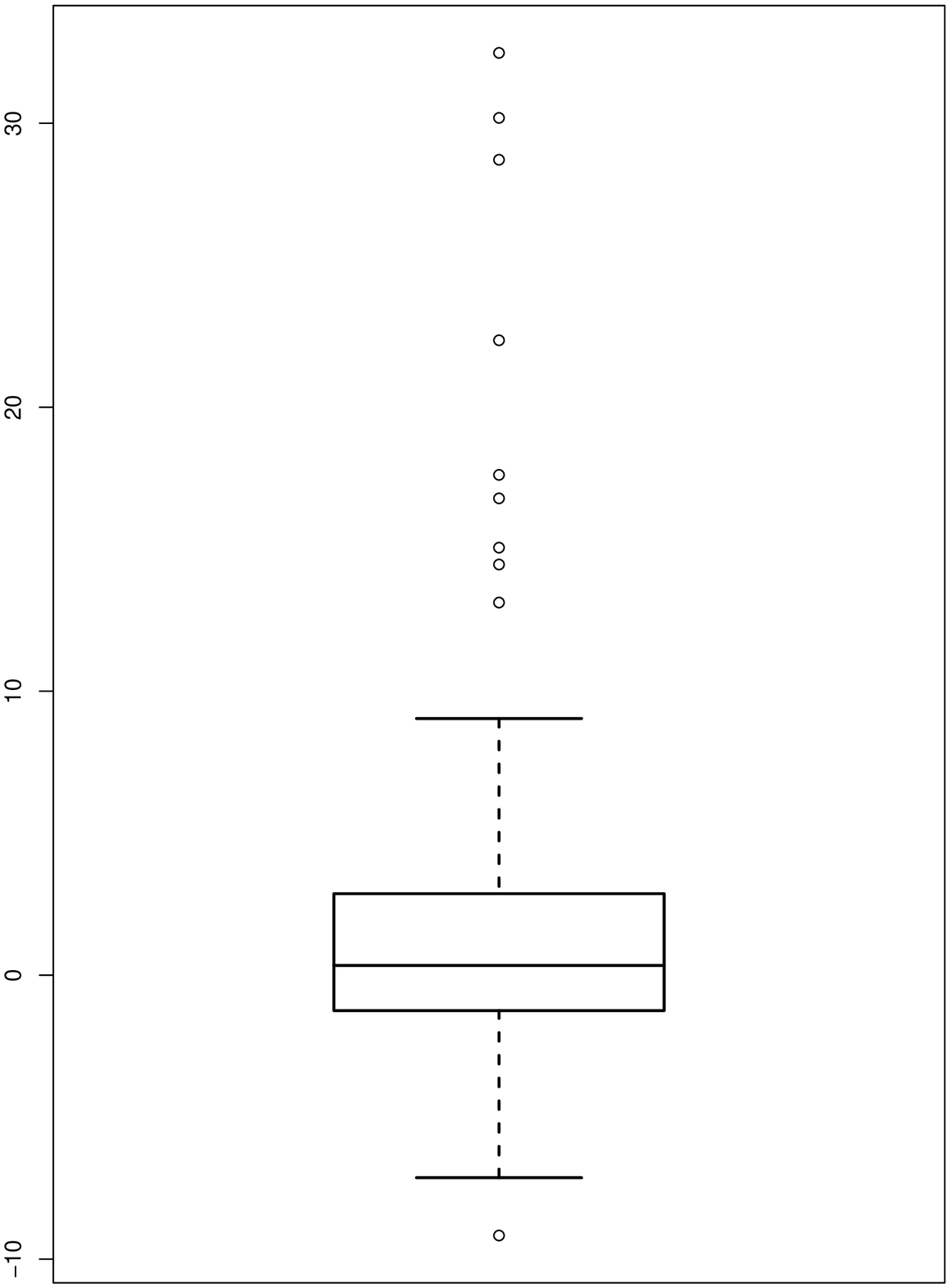} 
\end{tabular}
\end{center}
\caption{Comparison of hard and soft thresholding.  Left: 
$m(x)=5x_1^2 x_2^2$, $d=10$ and $\sigma = 0.5$; right: 
$m(x)=2(x_1+1)^3 + 2\sin(10\, x_2)$, $d=10$ and $\sigma = 1$.  The
hard and soft thresholding versions of the rodeo were compared on 100
randomly generated datasets, with a single random test point $x$ chosen
for each; $\beta = 0.9$.  The plots show two views of the difference of losses,
$\left(\tilde m_\text{hard}(x) - m(x)\right)^2 -
\left(\tilde m_\text{soft}(x) - m(x)\right)^2$; positive values
indicate an advantage for soft thresholding.}
\label{fig::hardsoft}
\end{figure}

We have taken the estimate 
\begin{equation}
\hat{D}_j(h) = Z_j(h)I( |Z_j(h)| > \lambda_j)
\end{equation}
with the result that
\begin{equation}
\tilde{m}(x) = \hat{m}_{\LL}(x) - \int_0^1 \langle \hat{D}(s), \overdot{h}(s)\rangle ds = \hat{m}_{h\star}(x).
\end{equation}
There are many possible generalizations.
First, we can replace $\hat{D}$ with the soft-thresholded estimate
\begin{equation}
\hat{D}_j(t)\, = \,{\rm sign}(Z_j(h))\,( |Z_j(h)| - \lambda_j)_+
\end{equation}
where the index $t$ denotes the $t^{\rm th}$ step of the
algorithm.  Since $h_j$ is updated multiplicatively as
$h_j \leftarrow \beta h_j$, the differential $dh_j(t)$ is given
by $dh_j(t) = (1-\beta)h_j$.  Using the resulting estimate of $D(t)$ and finite difference
approximation for $\overdot{h}(t)$ leads to the algorithm
detailed in Figure~\ref{fig::soft-rodeo}.

Figure~\ref{fig::hardsoft} shows a comparison of the hard and soft
thresholding versions of the rodeo on the example function
$m(x)=2(x_1+1)^3 + 2\sin(10\, x_2)$ in $d=10$ dimensions with $\sigma
= 1$; $\beta$ was set to $0.9$.  For each of 100 randomly generated datasets, a random test
point $x\sim \text{Uniform}(0,1)^{d}$ was generated, and the
difference in losses was computed:
\begin{equation}
\left(\tilde
m_\text{hard}(x) - m(x)\right)^2 - \left(\tilde m_\text{soft}(x) -
m(x)\right)^2\,.
\end{equation}
Thus, positive values indicate an advantage for soft thresholding,
which is seen to be slightly more robust on this example.

Another natural extension would be to consider more general paths than
paths that are restricted to be parallel to the axes.  We leave this
direction to future work.

\subsection{Global Rodeo}
\label{sec::global}

We have focused on estimation of $m$ locally at a point $x$.  The idea
can be extended to carry out global bandwidth and variable selection
by averaging over multiple evaluation points $x_1, \ldots, x_k$.
These could be points interest for estimation, could be randomly
chosen, or could be taken to be identical to the observed $X_i$s.

Averaging the $Z_j$s directly leads to a statistic whose mean for
relevant variables is asymptotically $k^{-1} h_j\sum_{i=1}^k
m_{jj}(x_i)$.  Because of sign changes in $m_{jj}(x)$, cancellations
can occur resulting in a small value of the statistic for
relevant variables.  To eliminate
the sign cancellation, we square the statistic.
Another way of deriving a global method would be to use
the statistic $\sup_x |Z_j^*(x)|$.

Let $x_1, \ldots, x_k$ denote the evaluation points.
Let
\begin{equation}
Z_j(x_i)= \sum_{s=1}^n Y_s G_j(X_s,x_i).
\end{equation}
Then define the statistic
\begin{equation}
T_j \equiv \frac{1}{k}\sum_{i=1}^k Z_j^2(x_i) = \frac{1}{k} Y^T P_j Y
\end{equation}
where $P_j = \GG_j \GG_j^T$, with 
$\GG_j(s,i) = G_j(X_s,x_i)$.

If $j\in R^c$ then we have $\E(Z_j(x_i))=o(1)$, so
it follows that, conditionally,
\begin{subequations}
\begin{eqnarray}
\E(T_j) &=&  \frac{\sigma^2}{k} \tr(P_j)+ o_P(1)\\
\V(T_j) &=&  \frac{2\sigma^4}{k^2} \tr(P_j P_j) + o_P(1).
\end{eqnarray}
\end{subequations}
We take the threshold to be
\begin{equation}
\lambda_j = \frac{\hat\sigma^2}{k} \tr(P_j) + 
2 \frac{\hat\sigma^2}{k} \sqrt{\tr(P_j P_j) \log(c_n n)} \;.
\end{equation}
Note that if $j > r$, we have
\begin{equation}
\E(T_j) = \frac{1}{k} \sum_i s_j^2(X_i) + O(h_0^6)
\end{equation}
but for $j\leq r$ we have
\begin{equation}
\E(T_j) = \frac{1}{k} \sum_i s_j^2(X_i) + O(h_0^2).
\end{equation}
We give an example of this algorithm in the following section,
leaving the detailed analysis of the asymptotics of this estimator
to future work.

\subsection{Greedier Rodeo and LARS}
\label{sec::lars}

The rodeo is related to least angle regression (LARS) \citep{lars:04}.
In forward stagewise linear regression, one performs variable
selection incrementally.  LARS gives a refinement where at each step
in the algorithm, one adds the covariate most correlated with the
residuals of the current fit, in small, incremental steps.  LARS takes
steps of a particular size: the smallest step that makes the largest
correlation equal to the next-largest correlation.  \cite{lars:04}
show that the lasso can be obtained by a simple modification of LARS.

The rodeo can be seen as a nonparametric version of forward stagewise
regression.  Note first that $Z_j$ is essentially the correlation
between the $Y_i$s and the $G_j(X_i,x,h)$s (the change in the
effective kernel).  Reducing the bandwidth is like adding in more of
that variable.  Suppose now that we make the following modifications
to the rodeo: (i) change the bandwidths one at a time, based on the
largest $Z^*_j = Z_j/\lambda_j$, (ii) reduce the bandwidth
continuously, rather than in discrete steps, until the largest $Z_j$
is equal to the next largest.  This version can then be thought of as
a nonparametric formulation of LARS.

In fact, we can go further and embed the rodeo within LARS to get a
fast nonparametric method.  We do this by replacing the derivatives of
the fit in the rodeo with differences.  Then we iterate variable
selection with bandwidth selection.

\begin{itemize}
\item Set $h=(\LL,\ldots, \LL)$.
Define $d$-dimensional pseudo-covariates
$\tilde{X}_i$, $i=1, \ldots, n$, by
\begin{equation}
\tilde{X}_i(j) = G_j(X_i,x,h),\;\; j=1, \ldots, d .
\end{equation}
Now run the LARS algorithm,
regressing the $Y_i$'s on the pseudo-covariates, up to some pre-defined stopping point.
This step essentially chooses relevant variables at the resolution of the
starting bandwidth $h=(\LL,\ldots, \LL)$.
\item
Define new pseudo-covariates
$\tilde{X}_i$, $i=1, \ldots, n$, by
\begin{equation}
\tilde{X}_i(j) = G_j(X_i,x,h') - G_j(X_i,x,h),\;\; j=1, \ldots, d
\end{equation}
where
$h'$ has $h_j$ replaced by $\beta h_j$.
Note that adding the $j^{\rm th}$ covariate corresponds to reducing the bandwidth from
$h_j$ to $\beta h_j$.
Now run the LARS algorithm up to some pre-defined stopping point.
\item Repeat the last step until a stopping criterion is satisfied.
\end{itemize}

One advantage of this method is that it can be implemented using
existing LARS software.  We leave the development of the theory for
this approach to future work.  Some examples of the greedy version
of this algorithm follow.

\subsubsection{Diabetes example}

\begin{figure}[t]
\begin{center}
\vskip-.4in
\begin{tabular}{c}
\includegraphics[width=3.2in,angle=-90]{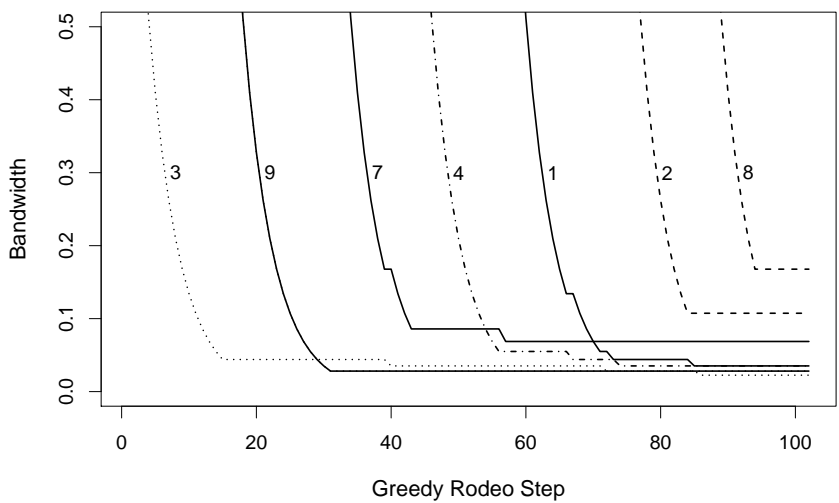} 
\end{tabular}
\end{center}
\vskip-10pt
\caption{Greedy rodeo on the diabetes data,
  used to illustrate LARS \citep{lars:04}.  A set of $k=100$ of the
  total $n=442$ points were sampled ($d=10$), and the bandwidth for
  the variable with largest average $|Z_j|/\lambda_j$ was reduced in
  each step.}
\label{fig::diabetes}
\end{figure}

Figure~\ref{fig::diabetes} shows the result of running the greedy
version of the rodeo on the diabetes dataset used by \cite{lars:04} to
illustrate LARS .  The algorithm averages $Z^*_j$ over a randomly
chosen set of $k=100$ data points, and reduces the bandwidth for
the variable with the largest value;  note that no estimate of
$\sigma$ is required.  The resulting variable ordering is
seen to be very similar to, but different from, the ordering obtained
from the parametric~LARS~fit.  The variables were selected in the
order 3 (body mass index), 9 (serum), 7 (serum), 4 (blood pressure), 1
(age), 2 (sex), 8 (serum), 5 (serum), 10 (serum), 6 (serum).  The LARS
algorithm adds variables in the order 3, 9, 4, 7, 2, 10, 5, 8, 6,
1. One notable difference is in the position of the age variable.

\subsubsection{Turlach's example}
\label{sec::turlach}

In the discussion to the LARS paper, Berwin Turlach \citep{Turlach}
gives an interesting example of where LARS and the lasso fails.  The
function is
\begin{equation} Y = \left(X_1-\frac{1}{2}\right)^2 + X_2 + X_3 + X_4 + X_5 + \varepsilon
\end{equation}
with ten variables $X_i \sim \text{Uniform}\,(0,1)$ and $\sigma=0.05$.
Although $X_1$ is a relevant variable, it is uncorrelated with $Y$,
and LARS and the lasso miss it.

\begin{figure}[t]
\begin{center}
\begin{tabular}{cc}
\hskip-60pt
\includegraphics[width=3.0in,angle=-90]{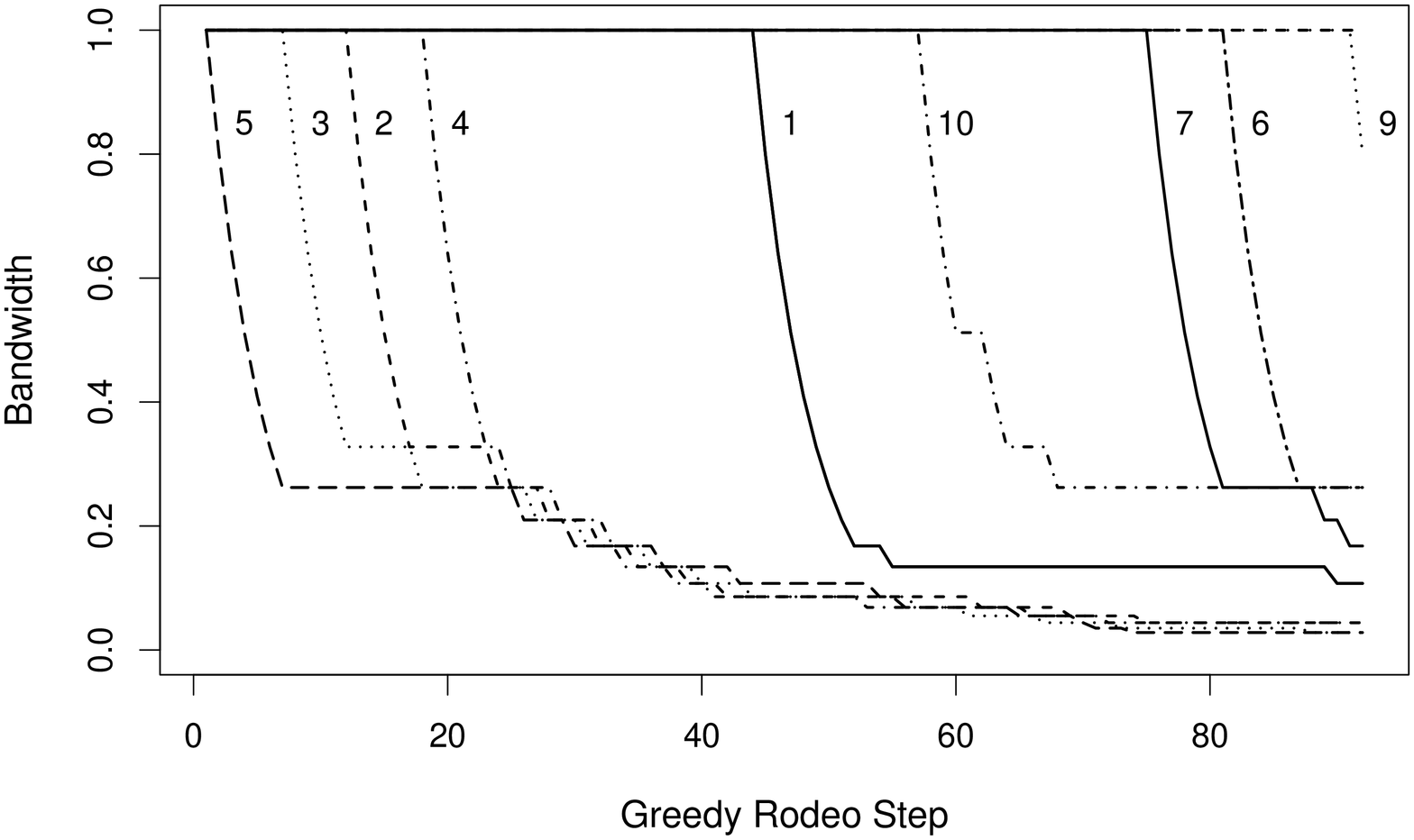} & 
\hskip-10pt
\includegraphics[width=1.4in,angle=-90]{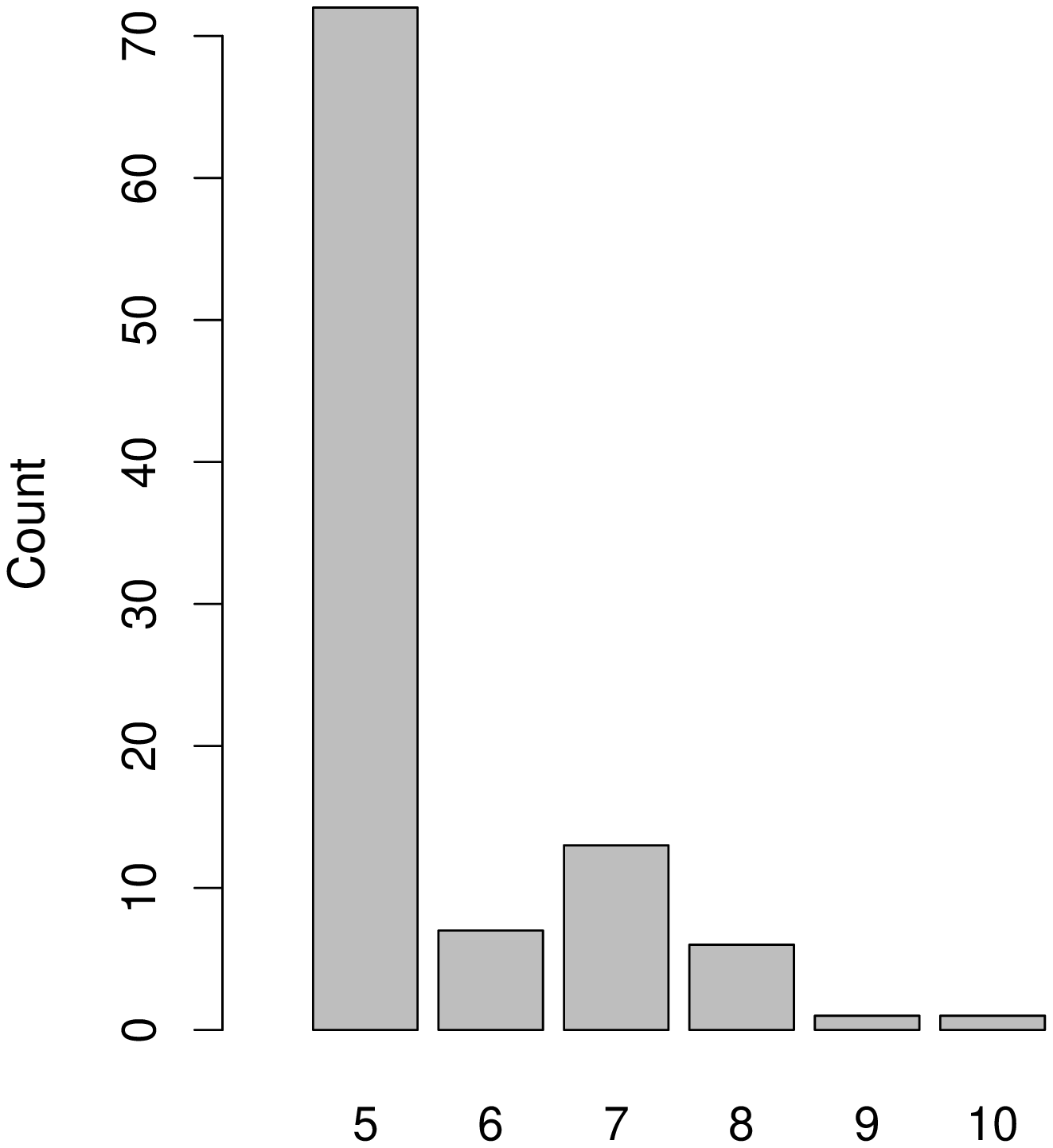}
\end{tabular}
\end{center}
\caption{Left: A typical run of the greedy algorithm on Turlach's
example.  The bandwidths are first reduced for variables $x_2, x_3,
x_4, x_5$, and then the relevant, but uncorrelated with $Y$ variable
$x_1$ is added to the model; the irrelevant variables enter the model
last.  Right: Histogram of the position at which variable $x_1$ is
selected, over 100 runs of the algorithm}
\label{fig::turlach}
\end{figure}

Figure~\ref{fig::turlach} shows the greedy algorithm on this example,
where bandwidth corresponding to the largest average $Z^*_j$ is
reduced in each step.  We use kernel regression rather than local
linear regression as the underlying estimator, without first
subtracting off a Lasso fit.  The variables $x_2$, $x_3$, $x_4$, $x_5$
are linear in the model, but are selected first in every run.
Variable $x_1$ is selected fifth in 72 of the 100 runs; a typical run
of the algorithm is shown in the left plot.  In contrast, as discussed
in \cite{Turlach}, LARS selects $x_1$ in position 5 about 25\% of the
time.

\begin{figure}[t]
\begin{center}
\begin{tabular}{cc}
\includegraphics[width=2.6in,angle=0]{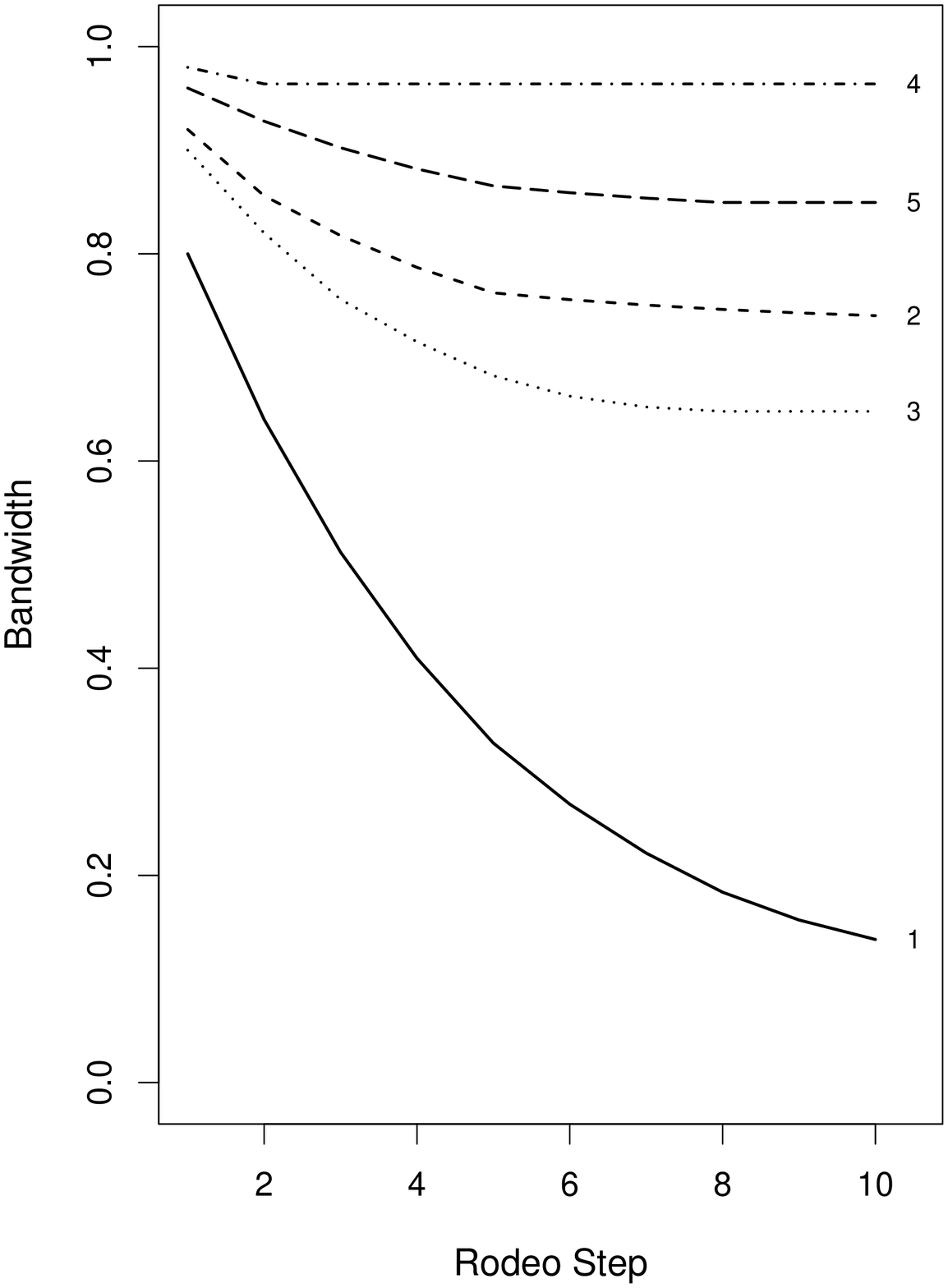} 
\includegraphics[width=2.6in,angle=0]{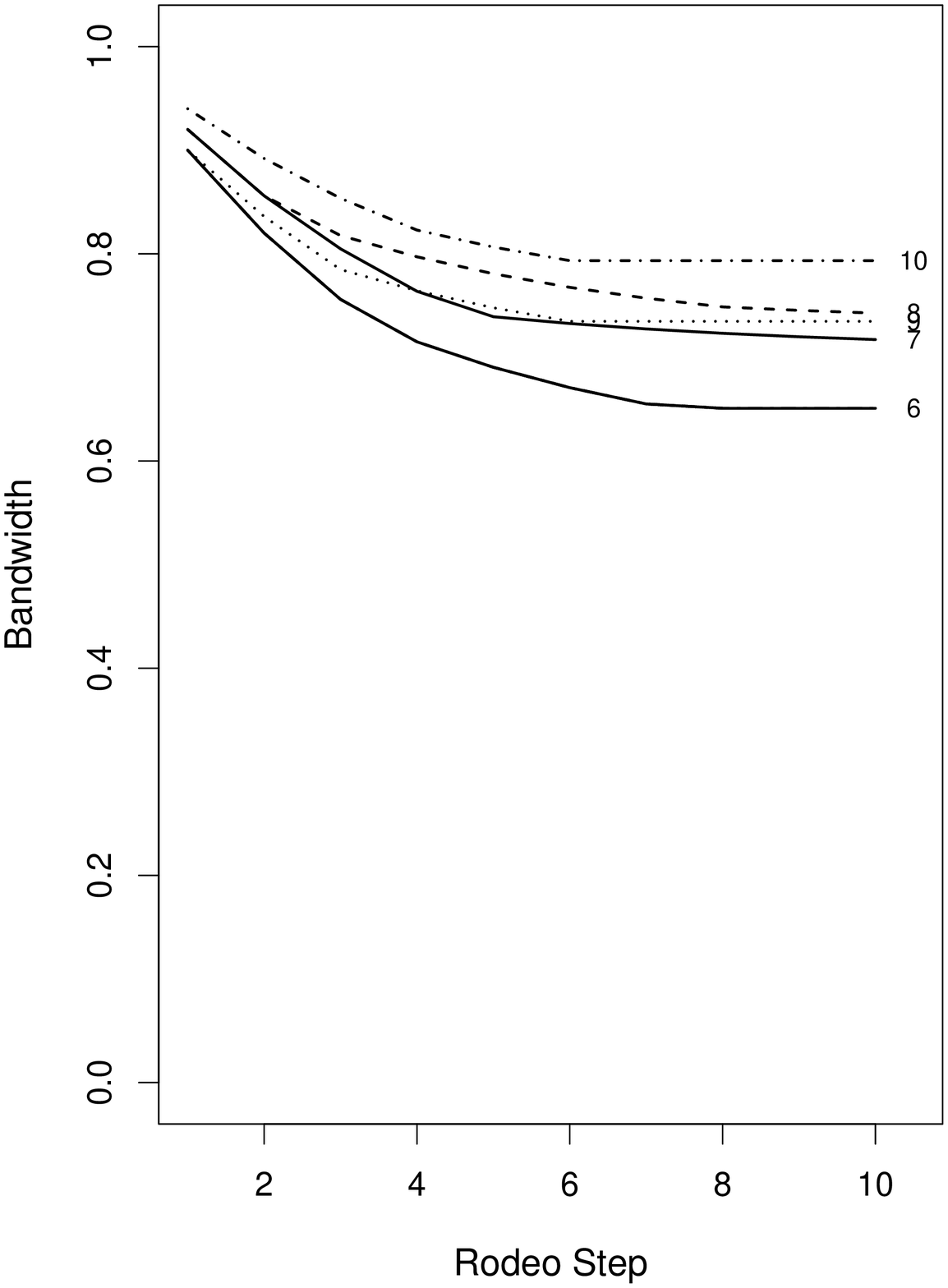} 
\end{tabular}
\end{center}
\caption{The global rodeo averaged over 10 runs on Turlach's example.
The left plot shows the bandwidths
for the five relevant variables. Since the linear effects
(variables two through five) have been subtracted off,
bandwidths $h_2,h_3,h_4,h_5$ are not shrunk.
The right plot shows the bandwidths for the other, irrelevant, variables.}
\label{fig::turlach-subtract}
\end{figure}

Figure \ref{fig::turlach-subtract} shows bandwidth traces for this
example using the global algorithm described in
Section~\ref{sec::global} with $k= 20$ evaluation points randomly
subselected from the data, and $\sigma$ taken to be known.  Before
starting the rodeo, we subtract off a linear least squares fit, and
run the rodeo on the residuals.  The first plot shows $h_1, \ldots,
h_5$.  The lowest line is $h_1$ which shrinks the most since $m$ is a
nonlinear function of $x_1$.  The other curves are the linear effects.
The right plot shows the traces for $h_6, \ldots, h_{10}$, the
bandwidths for the irrelevant variables.

\section{Proofs of Technical Results}
\label{sec::proofs}

In this section we give the proofs of the results
stated in Section~\ref{sec::properties}.

\subsection{Proof of Theorem~\ref{thm:important}}

We follow the setup of \cite{Ruppert:1994}
except for two differences:
the irrelevant variables
have different leading terms in the expansions than relevant variables
and $d$ is increasing; this requires a more detailed asymptotic
analysis than is carried out by \cite{Ruppert:1994}.

Let
$\H =
\left(\begin{array}{cc}
     \H_R & {0}\\
     {0} & {0}
     \end{array}\right)
$ denote the Hessian of $m(x)$,
let $h_j^{(t)}$ denote the $j^{\rm th}$ bandwidth at step $t$ and
denote the bandwidth matrix by $H^{(t)} = \text{diag}((h_1^{(t)})^2,\ldots,(h_d^{(t)})^2)$.
Similarly, let $H_R^{(t)} = \text{diag}((h_1^{(t)})^2,\ldots,(h_r^{(t)})^2)$.
Define
\begin{equation}
\mu_j^{(t)}=\frac{\partial}{\partial h_j} \E\left[\hat m_{H^{(t)}}(x)
  - m(x) \,|\, X_1,\ldots, X_d\right]
\end{equation}
In the following we sometimes supress the superscript $t$.

Let $\Dm$ be the gradient of $m$ at $x$, 
and let
\begin{equation}
Q = ( (X_1-x)^T{\cal H}(X_1-x), \ldots, (X_n-x)^T{\cal H}(X_n-x))^T .
\end{equation}
Note that
$\Dm$ and $Q$
are only functions of the relevant variables.
Then
\begin{equation}
m(X_i)= m(x) + (X_i-x)^T \Dm + \frac{1}{2} Q_i + R_i
\end{equation}
where, using multi-index notation,
\begin{equation}
R_i = \frac{1}{6} \sum_{|\alpha|=3} 
 (X_i-x)^\alpha \int_0^1 D^\alpha m\left((1-s) x + s X_i\right) \,ds =
 \frac{1}{6} \sum_{|\alpha|=3} (X_i-x)^\alpha D^{\alpha}m (\xi_{i,\alpha})
\end{equation}
where $\xi_{i,\alpha}$ is some point between $X_i$ and $x$.

So, with $M=(m(X_1),\ldots, m(X_n))^T$,
\begin{equation}
M = X_x 
\begin{pmatrix}
m(x) \\
\Dm
\end{pmatrix}
+ \frac{1}{2} Q + R
\end{equation}
where
$R = (R_1, \ldots, R_n)^T$.
Since
$S_x X_x ( m(x) , \Dm)^T  = m(x)$,
the bias 
$b(x) = \E(\hat{m}_H(x))-m(x)$
is given by
\begin{subequations}
\begin{eqnarray}
b(x) &=& e_1^T S_x M - m(x) = \frac{1}{2} S_x Q + S_x R \\
&=& \frac{1}{2} e_1^T (X^T_x W_x X_x)^{-1} X_x^T W_x Q +  e_1^T (X^T_x W_x X_x)^{-1} X_x^T W_x R\\
&=& \frac{1}{2} 
e_1^T \Upsilon_n^{-1} \Gamma_n +  e_1^T \Upsilon_n^{-1}\frac{1}{n} X_x^T W_x R
\end{eqnarray}
\end{subequations}
where
$\Upsilon_n = n^{-1} (X_x^T W_x X_x)$ and
$\Gamma_n = n^{-1} (X_x^T W_x Q)$.

{\it Analysis of $\Upsilon_n$.}
We write
\begin{subequations}
\begin{eqnarray}
\Upsilon_n &=&
\begin{pmatrix}
\frac{1}{n} \sum_{i=1}^n W_i & \frac{1}{n} \sum_{i=1}^n W_i (X_i-x)^T\\
\frac{1}{n} \sum_{i=1}^n W_i (X_i-x) & \frac{1}{n} \sum_{i=1}^n W_i (X_i-x) (X_i-x)^T
\end{pmatrix} =
\begin{pmatrix}
A_{11} & A_{12} \\
A_{21} & A_{22}
\end{pmatrix}
\end{eqnarray}
\end{subequations}
where 
\begin{equation}
W_i \;=\; \prod_{j=1}^d \frac{1}{h_j} K\left(\frac{x_j -
  X_{ij}}{h_j}\right) 
\;=\; \prod_{j=1}^d K_{h_j} (x_j - X_{ij})
\end{equation}
Now, the variance of $W_i$ can be bounded as
\begin{equation}
\Var(W_i) \leq E(W_i^2) =
 \prod_{j=1}^d \frac{1}{h_j} \E\left[
\prod_{j=1}^d \frac {1}{h_j} K^2 \left(
\frac{x_j - X_{ij}}{h_j}\right)\right] \leq
 C_K^d \prod_{j=1}^d \frac{1}{h_j}
\end{equation}
where $C_K=\sup_u K(u)$.
Therefore, 
\begin{equation}
\Var(A_{11}) \;\leq\;
\frac{C_K^d}{n} \prod_{j=1}^d \frac{1}{h_j} \leq 
\frac{1}{n}\frac{C_K^d }{\beta^{ds}h_0^d} \equiv \frac{\Delta_s}{n}.
\end{equation}
Also, $W_i \leq \Delta_s$.
Note that
$\beta^{T_n d} n h_0^d/C_k^d = \Omega\left(n^{1-\alpha-\log\log\log n
  / \log \log n}\right)$, where $\alpha$ is
the constant appearing in \eqref{eq::alpha}.
Hence, by Bernstein's inequality, for any $\delta > 0$,
\begin{subequations}
\begin{eqnarray}
\mathbb{P}( | A_{11} - \mathbb{E}(A_{11})| > \epsilon) & \leq &
2\exp\left\{ -\frac{1}{2} \left(\frac{n\epsilon^2}{\Delta + \Delta
  \epsilon/3}\right)\right\} \\
&=& 2\exp\left\{ -\frac{1}{2} 
\left(\frac{n\epsilon^2 \beta^{dt}h_0^d}{C_K^d (1 +
  \epsilon/3)}\right)\right\}\\
&\leq& 
2\exp\left\{ -\frac{1}{2} 
\left(\frac{n\epsilon^2 \beta^{dT_n}h_0^d}{C_K^d (1 +
  \epsilon/3)}\right)\right\} \\
&\leq& 
2\exp\left\{ -\frac{1}{2}  n^{(1-\delta)} \epsilon^2 \right\}.
\end{eqnarray}
\end{subequations}
Therefore, uniformly for all $t$ and $j$,
$A_{11}^{(t)} = \mathbb{E}(A_{11}^{(t)}) + O_P(1/n^{(1/2)-\delta})$
for any $\delta >0$.
Let $D_2$ denote the Hessian of $f$.
For some $\tilde{v}$ between $x$ and $v$,
\begin{subequations}
\begin{eqnarray}
\E(A_{11}) &=& 
\int \frac{1}{h_1 h_2 \cdots h_d} \, K(H^{-1/2}(x-u)) f(u)\, du \\
&=& 
\int K(v) f(x-\Hh v) \, dv \\
&=& f(x) - \int K(v) (\Hh v)^T \nabla f(x) \, dv + \frac{1}{2} \int
K(v) (\Hh v)^T D_2(\tilde{v}) \Hh v \, dv \nonumber
\\
\\
&=& f(x)  + \frac{1}{2} \int K(v) (\Hh v)^T D_2(\tilde{v}) \Hh v \, dv \\
&=& f(x) + r_1
\end{eqnarray}
\end{subequations}
where from \eqref{eq::d2}
\begin{equation}
\label{eq:r1}
r_1 = \frac{1}{2} \int K(v) (\Hh v)^T D_2(\tilde{v}) \Hh v \, dv =
 O\left(d \ h_0^2\ 
\max_{a,b}\sup_x \left| \frac{\partial^2 f(x)}{\partial x_a\partial
 x_b}\right| \right) = O\left((h_j^{(t)})^2 / d\right)
\end{equation}
Therefore we have, uniformly over $1\leq t\leq T_n$, 
\begin{equation}
A_{11} = f(x)  + O\left((h_j^{(t)})^2 / d\right).
\end{equation}

Next consider $A_{21}$.
By a similar argument as above,
uniformly for all $t$ and $j$,
$A_{21}^{(t)} = \mathbb{E}(A_{21}^{(t)}) + O_P(1/n^{(1/2)-\delta})$
for any $\delta >0$.
Also,
\begin{subequations}
\begin{eqnarray}
\E(A_{21}) &=& 
\int \frac{1}{h_1 h_2 \cdots h_d} \, K(H^{-1/2}(x-u)) (u-x)f(u)\, du \\
&=& \int K(v) \Hh v f(x+\Hh v) \, dv \\
&=& f(x) \underbrace{\int K(v) \Hh v \, dv}_{0} + \int K(v) \Hh v (\Hh
v)^T D\,  dv 
\nonumber \\
&& \quad + 
\frac{1}{2}\underbrace{\int K(v) \Hh v (\Hh v)^T D_2 \Hh v \,dv}_{0} \nonumber \\
&&\quad +
\frac{1}{6}\int K(v) \Hh v \sum_{abc}h_a h_b h_c v_a v_b v_c D_3(a,b,c) \,dv
\end{eqnarray}
\end{subequations}
where $D_2$ and $D_3$ are the arrays of second and third derivatives of $f$,
the latter being evaluated at some point between $x$ and $u$.
So,
\begin{subequations}
\begin{eqnarray}
\E(A_{21}) &=& 
\int K(v) \Hh v (\Hh v)^T D \,dv + 
\frac{1}{6}\int K(v) \Hh v \sum_{abc}h_a h_b h_c v_a v_b v_c
D_3(a,b,c)\,dv
\nonumber\\
\\
&=& \nu_2(K)H D + r_2
\end{eqnarray}
\end{subequations}
where
\begin{subequations}
\begin{eqnarray}
r_2 &=& \frac{1}{6}\int K(v) \Hh v \sum_{abc}h_a h_b h_c v_a v_b v_c
D_3(a,b,c) \,dv\\
&=& O\left( h_0^4 \ d^2 \ 
\max_{a,b,c,d}\sup_x \left| \frac{\partial^3 f(x)}
     {\partial x_a\partial x_b\partial x_c}\right| \vone
     \right) = O\left(((h_j^{(t)})^4/d) \, \vone\right)
\label{eq:r2}
\end{eqnarray}
\end{subequations}
using \eqref{eq::d3} and the definition of $\nu_2$ in equation~\eqref{eq:nudef}.
Therefore, 
$A_{21} = \nu_2(K)H D + O_P\left(((h_j^{(t)})^4/d) \vone\right)$.

Now we turn to $A_{22}$.
Again we have convergence to its mean and
\begin{subequations}
\begin{eqnarray}
\E(A_{22}) &=& \int K(v) \Hh v (\Hh v)^T f(x+\Hh v) \,dv \\
           &=& f(x)\int K(v) \Hh v (\Hh v)^T  \,dv + 
               \frac{1}{2} \int K(v) \Hh v (\Hh v)^T (\Hh v)^T D  \,dv
           \nonumber\\
&&\quad +   \frac{1}{6} \int K(v) \Hh v (\Hh v)^T (\Hh v)^T D_2 (\Hh v)\\
&=&
\nu_2 f(x) H + R_3
\end{eqnarray}
\end{subequations}
where $R_3$ is a matrix
whose off-diagonal elements are of order 
\begin{equation}
\label{eq:r3off}
O\left( h_0^4
\max_{a,b}\sup_x 
\left| \frac{\partial^2 f(x)}{\partial x_a\partial x_b}\right| \right) = O\left((h_j^{(t)})^4/d^2\right)
\end{equation}
and whose diagonal elements
are of order
\begin{equation}
\label{eq:r3}
O\left(d h_0^4
\max_{a,b}\sup_x \left| \frac{\partial^2 f(x)}{\partial x_a\partial x_b}\right| \right) = 
O\left((h_j^{(t)})^4/d\right).
\end{equation}
Thus,
\begin{equation}
\Upsilon_n = \left(\begin{array}{cc}
                    f(x) + r_1  & \nu_2 D^\t H + r_2^T\\
                    \nu_2 H^T D + r_2 & \nu_2 f(x) H + R_3
                     \end{array}\right).
\end{equation}
We can now write
\begin{equation}
\Upsilon_n = A + v v^T + E
\end{equation}
where
\begin{equation}
A = \begin{pmatrix}
f(x) - 1 & 0  \\
0        & \nu_2 f(x) H
\end{pmatrix},
\end{equation}
\begin{equation}
v^\top = (1,\  \nu_2 D^\t H)
\end{equation}
and
\begin{equation}
E =
\begin{pmatrix}
r_1 & r_2^T\\
r_2 & R_3 - \nu_2^2 H D D^T H
\end{pmatrix} .
\end{equation}
To invert the matrix we use the following result.
Let $B$ be an invertible matrix and let $\Delta$ be any matrix such that
$|\Delta| \leq \epsilon |B|$ (interpreted elementwise).
Then,
\begin{equation}
| (B + \Delta)^{-1} - B^{-1} | \leq \epsilon |B^{-1}|\ |B|\ |B^{-1}| + O(\epsilon^2)
\end{equation}
where the absolute values and inequalities are elementwise.
Thus,
\begin{equation}
| \Upsilon_n^{-1} - (A+v v^T)^{-1} | \leq 
\epsilon |(A+v v^T)^{-1}|\ |A+v v^T|\ |(A+v v^T)^{-1}| + O(\epsilon^2)
\end{equation}
where $\epsilon = \max_{ij} |E_{ij}|$.  It now follows from
\eqref{eq:r1}, \eqref{eq:r2}, \eqref{eq:r3off}, \eqref{eq:r3} and a
union bound that
\begin{equation}
\max_{j,k} | \Upsilon_n^{-1}(j,k) - (A+v v^T)^{-1}(j,k) | = o_P(h_0^2).
\end{equation}

To compute $(A+v v^T)^{-1}$
we use
the matrix inversion lemma (Woodbury formula)
\begin{equation}
(A + vv^\top)^{-1} = A^{-1} - A^{-1} v ( 1+v^\top A^{-1}  v)^{-1}
  v^\top A^{-1} \, .
\end{equation}
Now,
\begin{equation}
A^{-1} = \begin{pmatrix}
\frac{1}{f(x)-1} & 0 \\
0 & \frac{H^{-1}}{\nu_2 f(x)}
\end{pmatrix},
\end{equation}
and
\begin{equation}
1 + v^\top A^{-1} v \;=\; 
\frac{f(x)}{f(x)-1} + \frac{\nu_2 D^\top H D}{f(x)} =
\frac{f(x)}{f(x)-1}\left( 1 + O(d h_0^2 \max D_2)\right)
\end{equation}
so that
\begin{equation}
(1 + v^\top A^{-1} v)^{-1} \;=\; 
\frac{f(x)-1}{f(x)} ( 1 + o(1)).
\end{equation}
Also,
\begin{equation}
A^{-1} v v^\top A^{-1} \;=\; 
\begin{pmatrix}
\frac{1}{(f(x)-1)^2} & \frac{D^\top}{f(x)(f(x)-1)} \\
\frac{D}{f(x)(f(x)-1)}  & \frac{D D^T}{f(x)^2}
\end{pmatrix}.
\end{equation}
Finally, we have that
\begin{equation}
\Upsilon_n^{-1} = (A + vv^\top)^{-1} + o_P(1) =
\begin{pmatrix}
\frac{1}{f(x)} & -\frac{D^T}{f^2(x)}\\
 -\frac{D}{f^2(x)} & \frac{H^{-1}}{\nu_2 f(x)} - \frac{D D^T (f(x)-1)}{f^3(x)}
\end{pmatrix}
 +  o_P(h_0^2).
\end{equation}
This completes the analysis of $\Upsilon_n$.

{\it Analysis of $\Gamma_n = \frac{1}{n} X_x^T W_x Q$.}
We can write
\begin{equation}
\Gamma_n \equiv \frac{1}{n} X_x^T W_x Q =
\left(\begin{array}{l}
       \frac{1}{n}\sum_{i=1}^n W_i (X_i -x)^T {\cal H}(X_i -x)\\
       \frac{1}{n}\sum_{i=1}^n (X_i-x)W_i (X_i -x)^T {\cal H}(X_i -x)\\
     \end{array}
     \right) =
\left(\begin{array}{l}
       \gamma_1 \\ \gamma_2
     \end{array}
     \right).
\end{equation}
Now,
\begin{subequations}
\begin{eqnarray}
\mathbb{E}(\gamma_1) &=&
\int K(v) (\Hh v)^T {\cal H} (\Hh v) f(x + \Hh v) \,dv\\
&=& \nonumber
f(x) \int K(v) (\Hh v)^T {\cal H} (\Hh v) \,dv + 
\int K(v) (\Hh v)^T {\cal H} (\Hh v) D^T (\Hh v)\,dv \\
&&\ \ \ \ \ +
\frac{1}{2} \int K(v) (\Hh v)^T {\cal H} (\Hh v) (\Hh v)^T D_2^T (\Hh
v) \,dv\\
&=& \nu_2 f(x) \tr(H \H_R) + O(d r h_0^4 \max_{ab} D_2(a,b))\\
&=& \nu_2 f(x) \tr(H \H_R) + O((h_j^{(t)})^4/d).
\end{eqnarray}
\end{subequations}
The stochastic analysis of
$\gamma_1$ from its mean is similar to the analysis of $A_{12}$ and we have
$\gamma_1 = \nu_2 f(x) \tr(H \H_R) + o_P\left((h_j^{(t)})^4/d\right)$.

Next,
\begin{subequations}
\begin{eqnarray}
\mathbb{E}(\gamma_2) &=&
\int K(v) (\Hh v) (\Hh v)^T {\cal H} (\Hh v) f(x + \Hh v) \,dv\\
&=& f(x) \int (\Hh v)K(v) (\Hh v)^T {\cal H} (\Hh v) \,dv \nonumber \\
&& \quad +\;
\int K(v) (\Hh v) (\Hh v)^T {\cal H} (\Hh v) D^T (\Hh v) \,dv \nonumber \\
&&\quad +\;
\frac{1}{2} \int K(v)(\Hh v) (\Hh v)^T {\cal H} (\Hh v) (\Hh v)^T
D_2^T (\Hh v) \, dv
\\
&=& \int K(v) (\Hh v) (\Hh v)^T {\cal H} (\Hh v) D^T (\Hh v) \,dv \nonumber \\
&& \quad + \; \frac{1}{2} \int K(v)(\Hh v) (\Hh v)^T {\cal H} (\Hh v)
(\Hh v)^T D_2^T (\Hh v) \, dv
\\
&=& \Lambda + O(d r h_0^{5} \max_{ab} D_2(a,b)\, \vone)\\
&=& \Lambda + O\left((h_0^5/d)\vone\right)
\end{eqnarray}
\end{subequations}
where $\Lambda$ is a $d$-vector
with $j^{\rm th}$ component
\begin{equation}\label{eq::Lambda}
\Lambda_j = h_j^4 \nu_4 {\cal H}_{jj} f_j + I(j\leq r)\nu_2^2 h_j^2
\sum_{k=1}^r h_k^2 {\cal H}_{jk} f_k +
\nu_2^2 h_j^2 f_j \sum_{k=1}^r h_k^2 {\cal H}_{kk}.
\end{equation}
Thus,
\begin{equation}
\Gamma_n = \left(\begin{array}{c}
              \nu_2 f(x) \tr(H \H_R)\\
              \Lambda
             \end{array}
             \right) + O_P\left((h_0^5/d)\vone\right).
\end{equation}

{\it Analysis of remainder 
$e_1^T \Upsilon_n^{-1}\frac{1}{n} X_x^T W_x R$.}
We can write
\begin{equation}
\frac{1}{n} X_x^T W_x R \;=\;
\left(\begin{array}{l}
       \frac{1}{n}\sum_{i=1}^n W_i R_i\\
       \frac{1}{n}\sum_{i=1}^n (X_i-x)W_i R_i
     \end{array}
     \right) =
\left(\begin{array}{l}
       \delta_1 \\ \delta_2
     \end{array}
     \right).
\end{equation}
Let $m_3$ denote the array of third derivatives of $m$ and let
$\tilde{v}_j = h_j v_j$. Then,
\begin{subequations}
\begin{eqnarray}
\mathbb{E}(\delta_1) &=&
\int K(v) \, m_3(\xi) \sum_{1\leq a,b,c\leq
  r}\tilde{v}_a\tilde{v}_b\tilde{v}_c f(x + \tilde{v}) \, dv\\
&=&\nonumber
f(x)\int K(v) \, m_3(\xi) \sum_{1\leq a,b,c\leq
  r}\tilde{v}_a\tilde{v}_b\tilde{v}_c \, dv +
\int K(v) \, m_3(\xi) \sum_{1\leq a,b,c\leq
  r}\tilde{v}_a\tilde{v}_b\tilde{v}_c  \ (\tilde{v}^T D) \, dv \\
&& \\
&=& m_3(x) \int K(v) \sum_{1\leq a,b,c\leq
  r}\tilde{v}_a\tilde{v}_b\tilde{v}_c \ ( \tilde{v}^T D) \, dv + o\left((h_j^{(t)})^4\right)
\end{eqnarray}
\end{subequations}
and so
\begin{equation}
|\mathbb{E}(\delta_1)| =
O\Biggl(\sup_x |m_3(x)|\ \sup_x |D(x)|\ \sum_{j,k\leq r} h_j^2
h_k^2\Biggr) 
= O\Biggl(\sum_{j,k\leq r} h_j^2 h_k^2\Biggr) 
\end{equation}
Similarly,
\begin{subequations}
\begin{eqnarray}
\mathbb{E}(\delta_2) &=&
\int K(v) \, m_3(\xi) \, \tilde{v}\sum_{1\leq a,b,c\leq
  r}\tilde{v}_a\tilde{v}_b\tilde{v}_c f(x + \tilde{v}) \, dv\\
&=&
f(x) \int K(v) \, m_3(\xi)\, \tilde{v}\sum_{1\leq a,b,c\leq
  r}\tilde{v}_a\tilde{v}_b\tilde{v}_c \, dv + o\left((h_j^{(t)})^4\right)\\
&=& O\Biggl(\sup_x |m_3(x)| \sum_{j,k\leq r} h_j^2 h_k^2\Biggr) .
\end{eqnarray}
\end{subequations}
Hence,
\begin{equation}
e_1^T \Upsilon_n^{-1}\frac{1}{n} X_x^T W_x R =
O_P\Biggl(\sum_{j,k\leq r} h_j^2 h_k^2\Biggr) .
\end{equation}

Putting all of this together,
we conclude that the leading term of $b(x)$ is
\begin{equation}
\frac{1}{2} {\nu_2 \tr(H {\cal H}_R)} - \frac{1}{2}\frac{D^T\Lambda}{f^2(x)}
\end{equation}
where $\Lambda$ is given in (\ref{eq::Lambda}).
Taking the derivative with respect to bandwidth $h_j$ 
gives the statement of the theorem.

\begin{remark}
Special treatment is needed if $x$ is a boundary point;
see Theorem 2.2 of Ruppert and Wand (1994).
\end{remark}

\subsection{Proof of Theorem~\ref{thm:var}}

Let $\ell$ denote the first row of $S_x$.
Then, with $\xi \sim N(0,\sigma^2)$,
\begin{subequations}
\begin{eqnarray}
\hat{m}_H(x)    &=&  \sum_i \ell_i Y_i =  \sum_i \ell_i m(X_i) +  \sum_i \ell_i \epsilon_i\\
&\stackrel{d}{=} & \sum_i \ell_i m(X_i) + \xi\  \sqrt{\sum_i \ell_i^2}\\
&=& \sum_i \ell_i m(X_i) +  \frac{\Lambda}{\sqrt{n h_1 \cdots h_d}} \ \xi
\end{eqnarray}
\end{subequations}
where 
\begin{equation}
\Lambda = \sqrt{n h_1 \cdots h_d\sum_i \ell_i^2}.
\end{equation}
Thus,
\begin{equation}
{\rm Var}(Z_j(t)|X_1, \ldots, X_n) = 
\sigma^2 {\rm Var}\left( \frac{\partial}{\partial h_j} \frac{\Lambda}{\sqrt{n h_1 \cdots h_d}} \right).
\end{equation}
Now we find an asymptotic approximation for $\Lambda$.

Recall that
\begin{equation}
S_x =  \left(\frac{1}{n}X_x^T W_x X_x\right)^{-1} \frac{1}{n}X_x^T W_x 
\end{equation}
and from our previous calculations
\begin{equation}
\Upsilon_n^{-1} = \left(\frac{1}{n}X_x^T W_x X_x\right)^{-1} = 
\begin{pmatrix}
\frac{1}{f(x)} & -\frac{D^T}{f^2(x)}\\
 -\frac{D}{f^2(x)} & \frac{H^{-1}}{\nu_2 f(x)} - \frac{D D^T (f(x)-1)}{f^3(x)}
\end{pmatrix}
 + O_P\left((h_j^{(t)})^2/d\right).
\end{equation}
Note that
$\sum_i \ell_i^2$ is the (1,1) entry of
$S_x S_x^T$.
But
\begin{subequations}
\begin{eqnarray}
S_x S_x^T &=& (\Upsilon^{-1} \frac{1}{n} X_x^T W_x )(\Upsilon^{-1} \frac{1}{n} X_x^T W_x )^{T}\\
&=& \frac{1}{n^2} \Upsilon^{-1} X_x^T W_x^2 X_x \Upsilon^{-1}\\
&=& \frac{1}{n} \Upsilon^{-1} 
              \left(\begin{array}{cc}
               \frac{1}{n}\sum_i W_i^2 & \frac{1}{n}\sum_i (X_i-x)^T W_i^2\\
                \frac{1}{n}\sum_i (X_i-x) W_i^2 & \frac{1}{n}\sum_i (X_i-x)(X_i-x)^T  W_i^2
                \end{array}
                 \right)
                 \Upsilon^{-1} .
\end{eqnarray}
\end{subequations}
So $\Lambda^2$ is the $(1,1)$ entry of
\begin{subequations}
\begin{eqnarray}
&&  \Upsilon^{-1} 
              \left(\begin{array}{cc}
               \frac{h_1\cdots h_d}{n}\sum_i W_i^2 & \frac{h_1\cdots h_d}{n}\sum_i (X_i-x)^T W_i^2\\
                \frac{h_1\cdots h_d}{n}\sum_i (X_i-x) W_i^2 & \frac{h_1\cdots h_d}{n}\sum_i (X_i-x)(X_i-x)^T  W_i^2
                \end{array}
                 \right)
                 \Upsilon^{-1}\\
  &=& \Upsilon^{-1} 
              \left(\begin{array}{cc}
               a_{11} & a_{12}\\
               a_{21} & a_{22}
                \end{array}
                 \right)
                 \Upsilon^{-1}
\end{eqnarray}
\end{subequations}
\def\u{\bar\nu}
Next, as in our earlier analysis
\begin{subequations}
\begin{eqnarray}
a_{11} &=& \int  K^2(v) \, f(x - \Hh v)\, dv\\
&=& f(x) \int K^2(v) dv + O_P(h_0^2/d).
\end{eqnarray}
\end{subequations}
Similarly,
\begin{equation}
a_{21} = \u_2 H D + O_P\left(((h_j^{(t)})^4/d) \vone\right)
\end{equation}
where $\u_2 I = \int v v^T K^2(v) \, dv$ and
\begin{equation}
a_{22} = f(x) \u_2 H + E
\end{equation}
where $\max_{ij} E_{ij} =  O_P\left((h_j^{(t)})^4\right)$.
Hence, the leading order expansion of $\Lambda^2$ is given by
\begin{equation}
\frac{\int K^2(v) dv}{f(x)} - \u_2\frac{D^T H D}{f^3(x)} + o(\tr(H)) = 
\frac{\int K^2(v) dv}{f(x)} + O(\tr(H))
\end{equation}
Taking the derivative with respect to $h_j$ we thus conclude
that
\begin{equation}
{\rm Var}(Z_j(t)|X_1, \ldots, X_n) = 
\frac{\sigma^2 \int K^2(v)\, dv}{f(x) h_j^2}\frac{1}{n h_1\cdots h_d}(1 + o_P(1)).
\end{equation}

\subsection{Proof of Theorem~\ref{thm}}

This result characterizes the asymptotic running time, 
selected bandwidths, and risk of the algorithm.  We restate 
the assumptions made on 
the functions $m$ and $f$ in Section~\ref{sec::properties}.

\begin{enumerate}
\item[(A1)]
For each $j > r$, 
\begin{equation}
\nabla_j \log f(x) = O\left( \frac{1}{n^{1/4}}\right)
\end{equation}
\item[(A2)]
\begin{equation}
\liminf_n \min_{1\leq j\leq r} |m_{jj}(x)| > 0 \,.
\end{equation} 
\end{enumerate}

To prove the theorem we make use of a version of Mill's
inequality, modified for non-zero mean random variables.

%
%
%
%

\begin{lemma}\label{lemma::mills}
Let $Z\sim N(\mu,\sigma^2)$.
If $\lambda \geq 2\mu$
and $\lambda^2 \geq 2 \sigma^2$ then
\begin{eqnarray}
\mathbb{P}(|Z| > \lambda)  \leq 
\frac{5\lambda}{\sigma}\exp\left\{ - \frac{\lambda^2}{8 \sigma^2}\right\} .
\end{eqnarray}
Moreover, if $\lambda \geq 5 \sigma$ then
\begin{eqnarray}
\mathbb{P}(|Z| > \lambda)  \leq 
\exp\left\{ - \frac{\lambda^2}{16 \sigma^2}\right\} .
\end{eqnarray}
\end{lemma}

\begin{proof}
Without loss of generality, assume $\mu >0$.
Then,
\begin{subequations}
\begin{eqnarray}
\mathbb{P}(|Z| > \lambda) & \leq & 2\mathbb{P}(Z > \lambda) = 2\mathbb{P}\left(\frac{Z-\mu}{\sigma} > \frac{\lambda - \mu}{\sigma}\right)\\
&\leq & \frac{2\sigma}{\lambda - \mu} \exp\left\{ - \frac{(\lambda - \mu)^2}{2 \sigma^2}\right\}  \equiv  B(\mu).
\end{eqnarray}
\end{subequations}
Now $B(\mu) = B(0) + \mu B'(\tilde{\mu})$ for some
$0\leq \tilde{\mu} \leq \mu$ and
\begin{eqnarray}
B'(\mu) = \frac{2\sigma}{\lambda - \mu} \exp\left\{ - \frac{(\lambda - \mu)^2}{2 \sigma^2}\right\}  
\left(\frac{\lambda - \mu}{\sigma^2} + \frac{1}{\lambda - \mu}\right).
\end{eqnarray}
Hence,
\begin{eqnarray}
B'(\tilde\mu) \leq \frac{2\sigma}{\lambda - \mu} \exp\left\{ - \frac{(\lambda - \mu)^2}{2 \sigma^2}\right\}  
\left(\frac{\lambda }{\sigma^2} + \frac{1}{\lambda - \mu}\right).
\end{eqnarray}
When $\lambda \geq 2\mu$,
$1/(\lambda - \mu) \leq 2/\lambda$ and
$(\lambda - \mu)^2 \geq \lambda^2/4$ so that if $\lambda^2 \geq 2\sigma^2$ then
\begin{eqnarray}
B'(\tilde\mu) \leq \frac{4\sigma}{\lambda} \exp\left\{ - \frac{\lambda^2}{8 \sigma^2}\right\}  
\left(\frac{\lambda }{\sigma^2} + \frac{2}{\lambda}\right) \leq
 \frac{8}{\sigma} \exp\left\{ - \frac{\lambda^2}{8 \sigma^2}\right\}.  
\end{eqnarray}
Thus,
\begin{eqnarray}
\mathbb{P}(|Z| > \lambda) & \leq & 
\frac{2\sigma}{\lambda}\exp\left\{ - \frac{\lambda^2}{2 \sigma^2}\right\}   + 
 \frac{8\mu}{\sigma} \exp\left\{ - \frac{\lambda^2}{8 \sigma^2}\right\} \leq
\frac{5\lambda}{\sigma}\exp\left\{ - \frac{\lambda^2}{8 \sigma^2}\right\} .
\end{eqnarray}
The last statement follows since
$5 x  e^{-x^2/8} \leq e^{-x^2/16}$ for all $x \geq 5$. 
\end{proof}

\begin{proofof}{Theorem~\ref{thm}}
First consider $j > r$.  
Let $V_t=\{j>r:\ h_j = \LL \beta^t \}$ be the set of
irrelevant dimensions that are active at stage $t > 1$
of the algorithm.
From Assumption (A1) it follows that, for sufficiently large $n$,
\begin{equation}
\lambda_j \geq | L_j^{(t)} | = O(h_0^3 \, n^{-1/2})
\end{equation}
for $j > r$.
Also, $\lambda_j\geq 5s_j^{(t)}$.  From Theorems 5.1 and 5.2,
we have that, with probability tending to one uniformly over
$j > r$, $\lambda_j\geq 2 |\mu_j^{(t)}| $ and
$\lambda_j^2\geq 2v_j^{(t)}$, and $v_j^{(t)} / (s_j^{(t)})^2 = 1 + o(1)$.
Thus, by Lemma \ref{lemma::mills},
\begin{subequations}
\begin{eqnarray}
\P( |Z_j| > \lambda_j,\ {\rm for\ some\ }j \in V_t) & \leq &
\sum_{j\in V_t} \P( |Z_j| > \lambda_j) + o(1) \\
& \leq & d \exp\left(-\lambda_j^2/(16 v_j)\right) + o(1) \\
&=& d \exp\left(-\lambda_j^2(1+o(1))/(16 s_j^2)\right) + o(1)  \to 0.
\end{eqnarray}
\end{subequations}
Therefore, with probability tending to one,
$h_j = h_0$  for each $j>r$, meaning that the bandwidth for each irrelevant
dimension is frozen in the first step in the algorithm.

Now consider  $j\leq r$.
By (\ref{eq::A}), and Theorem \ref{thm:important}, for all large $n$,
\begin{equation}
|\mu_j| \geq |L_j^{(t)}| - |\mu_j - L_j^{(t)}| = |L_j^{(t)}| - O_P((h_j^{(t)})^4/d) \geq
c h_j^{(t)} |m_{jj}(x)|
\end{equation}
for some $c>0$.
Without loss of generality, assume that
$c h_j m_{jj}(x) >0$.
We claim that in iteration $t$ of the algorithm, if 
\begin{equation}
\label{eq:tcond}
t \; \leq \; \frac{1}{4+r} 
\log_{1/\beta} \left(\frac{c^2 n \Amin^2 \LL^{d+4}}{8C  \log(n
  c_n)}\right)
= t_0
\end{equation}
then 
\begin{equation}
\P(h_j = \LL\beta^t,\;\mbox{for all $j\leq r$}) \longrightarrow 1 .
\end{equation}

To show this, first note that \eqref{eq:tcond} can be written as
\begin{equation}
\left(\frac{1}{\beta}\right)^ {t(4+r)} 
\;\leq \; \frac{c^2 n \Amin^2 \LL^{d+4}}{8C \log (c_n  n)} \,.
\end{equation}
Except on an event of vanishing probability, we have shown above that
\begin{equation}
\prod_{j > r} \frac{1}{h_j} \;=\; \frac{1}{\LL^{(d-r)}}\,.
\end{equation}
So on the complement of this event, if each relevant dimension
is active at step $s\leq t$, we have
\begin{subequations}
\begin{eqnarray}
\frac{\lambda_j^2}{h_j^2} & = &\frac{2 s_j^2 \log (n c_n)}{h_j^2} =
\frac{2 C\log(n c_n)}{n h_j^4}  \prod_i \frac{1}{h_i}\\
& \leq & \frac{2 C \log(n c_n) }{n\LL^{d+4}}
\left(\frac{1}{\beta}\right)^{(4+r)t}  \leq  \frac{ c^2  \Amin^2}{4}  \\
&\leq&  \frac{ c^2 m_{jj}(x)^2}{4} 
\end{eqnarray}
\end{subequations}
which implies that
\begin{equation}
c m_{jj}(x) h_j \geq 2\lambda_j 
\end{equation}
and hence
\begin{equation}\label{eq:cb}
\frac{c m_{jj}(x) h_j - \lambda_j}{s_j} \geq \frac{\lambda_j}{s_j} = 
\sqrt{2 \log (n c_n)}
\end{equation}
for each $j\leq r$.
Now, 
\begin{subequations}
\begin{eqnarray}
\P(\text{rodeo halts}) &=&
\P(|Z_j| < \lambda_j\;\mbox{for all $j\leq r$}) \leq
\P(|Z_j| < \lambda_j\;\mbox{for some $j\leq r$})\\
&\leq&\sum_{j\leq r} \P(|Z_j| < \lambda_j)  \leq \sum_{j\leq r}\P(Z_j < \lambda_j) \\
&\leq& \sum_{j\leq r}\P\left(\frac{Z_j-\mu_j}{s_j} > \frac{\mu_j -   \lambda_j}{s_j}\right)\\
&\leq& \sum_{j\leq r}\P\left(\frac{Z_j-\mu_j}{s_j} > \frac{c m_{jj}(x) h_j -   \lambda_j}{s_j}\right)\\
&\leq & \frac{r}{n c_n \sqrt{2 \log (n c_n)}}  \,.
\end{eqnarray}
\end{subequations}
Finally, summing
over all iterations $s \leq t$ gives
\begin{subequations}
\begin{eqnarray}
\P\left(\bigcup_{s\leq t} \bigcup_{j\leq r} 
\left\{|Z^{(t)}_j| < \lambda^{(t)}_j\right\}\right)  
& \leq & \frac{tr}{n c_n \sqrt{2 \log (n c_n)}} \\
& \leq & 
\frac{\displaystyle \frac{r}{r+4} \log_{1/\beta}
\left(\frac{c^2 n \Amin^2 \LL^{d+4}}
{8 C \log(n c_n )}\right)}{n c_n \sqrt{2 \log(n c_n)}} \longrightarrow
0.
\end{eqnarray}
\end{subequations}
Thus, the bandwidths $h_j$ for
$j \leq r$ satisfy, with high probability,
\begin{subequations}
\begin{eqnarray}
h_j \;=\; \LL\beta^{t} &\leq& \LL\beta^{t_0} \\
&=&\LL \left(\frac{8C \log(n c_n)}{c^2
  \Amin^2 n \LL^{d+4}}\right)^{1/(4+r)} \\
&=& n^{-1/(4+r)} \left(\frac{8C \log(n c_n)}{c^2  \Amin^2 \LL^{d-r}}\right)^{1/(4+r)} \,.
\end{eqnarray}
\end{subequations}
Now, since $d\log d = O(\log n)$ it follows that
$d = O(\log n / \log \log n)$ and 
\begin{subequations}
\begin{eqnarray}
\frac{1}{\LL^d} &=& O\left((\log\log n)^d\right) \\
&=& O\left( n^{\textstyle \frac{\log\log\log n}{\log\log n}}\right) \\
&=& O\left(n^\epsilon\right) 
\end{eqnarray}
\end{subequations}
for any $\epsilon > 0$.
In particular, with probability approaching one, the algorithm
runs for a number of iterations $T_n$ bounded from below by
\begin{equation}
T_n \geq \frac{1}{4+r} \log_{1/\beta} (n^{1-\epsilon} a_n)
\end{equation}
where
\begin{equation}
a_n \;=\; \frac{c^2 \Amin^2}{8C \log(n c_n)} \;=\;
    \tilde \Omega(1)\,.
\end{equation}

We next show that the algorithm is unlikely to reach iteration
$s$, if 
\begin{equation}
\label{eq:scond}
s \;\geq\; \frac{1}{4+r} \log_{1/\beta} \left(n^{1+\delta} b_n\right)
\end{equation}
for any $\delta > 0$,
where $b_n = \tilde O(1)$ is defined below.  
From the argument
above, we know that except on an event of vanishing probability,
each relevant dimension $j\leq r$ has bandwidth no larger than
\begin{subequations}
\begin{eqnarray}
h_j &\leq& h_0 \beta^{\textstyle\,\frac{1}{4+r} \log_{1/\beta} (n^{1-\epsilon} a_n)} \\
   &=& \frac{h_0}{(n^{1-\epsilon} a_n)^{1/(4+r)}} .
\end{eqnarray}
\end{subequations}

Thus, if relevant dimension $j$ has bandwidth $h_j \leq h_0\beta^s$,
then from Theorem~\ref{thm:var} we have that
\begin{subequations}
\begin{eqnarray}
\frac{s_j^2}{A_j^2 h_j^2} &\geq& \frac{C}{A_j^2 n h_0^4 \beta^{4s}} 
 \frac{n^{r(1-\epsilon)/(4+r)}a_n^{r/(4+r)}}{h_0^r} \frac{1}{h_0^{d-r}} \\
&=& \frac{C}{A_j^2 n^{4/(4+r)+\delta}} \frac{a_n^{r/(4+r)}}{h_0^{4+d}}  \frac{1}{\beta^{4s}} \\
&\geq& \frac{C}{\Amax^2 n^{4/(4+r)+\delta}} \frac{a_n^{r/(4+r)}}{h_0^{4+d}}
 \frac{1}{\beta^{4s}} 
\end{eqnarray}
\end{subequations}
where $\delta = r\epsilon/(4+r)$.
Therefore, 
\begin{equation}
\frac{s_j^2}{A_j^2 h_j^2} \geq \log\log n
\end{equation}
in case
\begin{subequations}
\begin{eqnarray}
\left(\frac{1}{\beta}\right)^s &\geq& 
(n^{1+\delta} b_n)^{1/(4+r)} \\
&=&
n^{1/(4+r)+\delta'}  \left(\frac{\Amax^2 \log\log n}{C
  a_n^{r/(4+r)}}\right)^{1/4} \\
&\geq&
n^{1/(4+r)+\delta'}  \left(\frac{\Amax^2 h_0^{4+d} \log\log n}{C
  a_n^{r/(4+r)}}\right)^{1/4} 
\end{eqnarray}
\end{subequations}
which defines $b_n = \tilde O(1)$.
It follows that in iteration $s \geq \frac{1}{4+r} \log_{1/\beta}
\left(n^{1+\delta} b_n\right)$, the probability of a relevant variable having
estimated derivative $Z_j$ above threshold is bounded by
\begin{subequations}
\begin{eqnarray}
\P( |Z_j| > \lambda_j,\ {\rm for\ some\ }j \leq r) & \leq &
\sum_{j \leq r} \P\left( \frac{|Z_j|}{s_j} > \frac{\lambda_j}{s_j}\right)\\
&\leq& \sum_{j\leq r}
\Biggl(\frac{s_j}{\lambda_j}e^{-\lambda_j^2/(2s_j^2)} +   
 \frac{1}{4}\frac{\mu_j^2}{s_j^2}\Biggr)\\
&\leq & \frac{r}{n c_n \sqrt{2 \log (n c_n)}} + \frac{1}{4}\sum_{j\leq r}
 \frac{\mu_j^2}{s_j^2} \\
&\leq & \frac{r}{nc_n \sqrt{2 \log (n c_n)}} + \frac{r}{4\log\log n} \\
&=& O\left(\frac{1}{\log\log n}\right)
\end{eqnarray}
\end{subequations}
which gives the statement of the theorem.
\end{proofof}

\begin{proofof}{Corollary~\ref{cor}}
Given the bandwidths in \eqref{eq:bw1} and \eqref{eq:bw2}
we have that the squared (conditional) bias is given by
\begin{subequations}
\begin{eqnarray}
\mbox{Bias}^2(\hat{m}_{h^\star}) &=& \left(\textstyle
\sum_{j\leq r} A_j {h_j^\star}^2\right)^2
+ o_P(\text{tr}({H^\star}))\\
&=& \sum_{i,j\leq r} A_i A_j {h_i^\star}^2 {h_j^\star}^2 
+ o_P(\text{tr}({H^\star}))\\
&=& {O}_P(n^{-4/(4+r)+\epsilon})
\end{eqnarray}
\end{subequations}
by Theorem~\ref{thm}.  Similarly, 
from the proof of Theorem~\ref{thm:var} the (conditional) variance is 
\begin{subequations}
\begin{eqnarray}
\Var(\hat{m}_{h^\star}) &=& \frac{1}{n} \left(\prod_i \frac{1}{{h_i^\star}}\right)
\frac{R(K)}{f(x)} \sigma^2 (1+ o_P(1)) \\
&=& {O}_P(n^{-1+r/(r+4)+\epsilon}) \\
&=& {O}_P(n^{-4/(4+r)+\epsilon})
\end{eqnarray}
\end{subequations}
where $R(K) = \int K(u)^2 \, du$.  The result follows from
the bias-variance decomposition.
\end{proofof}

\subsection{Proof of Theorem~\ref{thm:noassumption}}

This result shows that assumption (A1) can be dropped if we modify the
algorithm by simply increasing the threshold slightly.
Examples of choices for $\rho_n$ are
$\rho_n = \delta h_0^2$ or
$\rho_n =  h_0^{3-\delta}$
for some small $\delta >0$.

Let $j>r$.  Without loss of generality assume that $\mu_j \geq 0$.
Recall that, if variable $j$ is still active in round $t$,
$\mu_j = ch_j^3 + O_P((h_j^{(t)})^2/d)$
for some constant~$c > 0$.
It follows from (\ref{eq::rho}) that, for all large $n$,
\begin{equation}
\lambda_j > \mu_j + s_j\sqrt{2\log(n c_n)}
\end{equation}
and hence
\begin{equation}
\frac{ s_j}{\lambda_j-\mu_j} e^{- (\lambda_j - \mu_j)^2/(2s_j^2)} \leq
\frac{1}{n c_n \sqrt{2\log(n c_n)}}.
\end{equation}
Then, with $\xi \sim N(0,1)$, we have
\begin{subequations}
\begin{eqnarray}
\mathbb{P}(|Z_j| > \lambda_j, \ {\rm for\ some\ j>r}) & \leq &
\sum_j  2\mathbb{P}\left( \xi > \frac{\lambda_j -\mu_j}{s_j}\right)\\
& \leq & \sum_j \frac{2s_j}{\lambda_j-\mu_j} e^{- (\lambda_j - \mu_j)^2/(2s_j^2)}\\
& \leq & \frac{2d}{n c_n \sqrt{2\log(n c_n)}} \longrightarrow 0.
\end{eqnarray}
\end{subequations}

Now consider $j \leq r$.
It follows from (\ref{eq::rho}) that, for all large $n$,
\begin{equation}
\frac{\mu_j - \rho_n\beta^{3t}}{s_j} > c \frac{\mu_j}{s_j}
\end{equation}
for some $c>0$.
Hence,
\begin{subequations}
\begin{eqnarray}
\mathbb{P}(Z_j< \lambda_j)&=&\mathbb{P}\left(\frac{Z_j-\mu_j}{s_j}< \frac{\lambda_j-\mu_j}{s_j}\right)\\
&=&\mathbb{P}\left(\xi > \frac{\mu_j - \lambda_j}{s_j}\right)\\
& \leq & \mathbb{P}\left(\xi > \frac{c\mu_j - \lambda_j'}{s_j}\right)
\end{eqnarray}
\end{subequations}
where
$\lambda_j' = s_j\sqrt{2\log(n c_n)}.$
The calculation then proceeds as before.

\section{Acknowledgments}

We thank Gilles Blanchard and the editors and referees for many helpful comments
on this work.

\vskip20pt
\bibliography{rodeo}

\end{document}